\documentclass[10pt,a4paper]{article}
\usepackage[utf8]{inputenc}
\usepackage{amsmath}
\usepackage{amsfonts}
\usepackage{amssymb}
\usepackage{url}
\usepackage{graphicx}
\usepackage{amsmath}
\usepackage{graphicx}
\usepackage{amsfonts}
\usepackage{amsthm}
\graphicspath{ {images/} }

\newtheorem*{theorem1}{\textit{Theorem 1}$^\prime$}

\newtheorem*{theorem3}{\textit{Theorem 3}$^\prime$}

\DeclareMathOperator{\Ima}{Im}

\author{Alexander Kushkuley \\ (Salesforce, akushkuley@salesforce.com, \\ kushkuley@gmail.com)       
}

\title{Block Approximation of Tall Sparse Matrices and Block-Givens Rotations}

\begin{document}
\maketitle

\begin{abstract}
\noindent Estimation of top singular values  is one of the widely used techniques and one of the intensively researched problems in Numerical Linear Algebra and  Data Science. We consider here two general questions related to this problem:
\begin{itemize}
	\item[] How top singular values are affected by zeroing out a sparse rectangular block of a matrix?
	\item[]  How much top singular values differ from top column norms of a tall sparse non-negative matrix ?
\end{itemize}   
\end{abstract}
\paragraph{}

\newtheorem{lemma}{Lemma}
\newtheorem{remark}{Remark}
\newtheorem{definition}{Definition}
\newtheorem{theorem}{Theorem}

\newtheorem{algo}{Algorithm}
\newtheorem{corollary}{Corollary}
\newtheorem{example}{Example}
\newcommand{\reals}{\mathbb{R}}
\newcommand{\complex}{\mathbb{C}}

\section{Introduction}
For a real matrix $ X $ of rank $ r $ its singular values are indexed in decreasing order, i. e.
it is always assumed that
\begin{equation}
 \sigma_1 \geq \sigma_2 \geq \cdots \geq \sigma_r > 0 \nonumber  
\end{equation}
Everywhere below $ \|X\| $ means an operator norm of a matrix $ X $ and by
$ \sigma_i(X) $ we denote  the $i$-th singular value of $ X $  so that $  \| X \| \equiv \sigma_1( X ) .$
It is  well known (cf. e.g. \cite{Ho},\cite{Go}) that the optimal   rank $k$ approximation $ X_k $ of the a matrix  $ X $  is characterized by the property 
\begin{equation}
\| X - X_k \| = \sigma_{k+1}( X )  \nonumber
\end{equation} 
Let 
\begin{equation} \label{R}
R = 
\left( \begin{array}{cc}
A & B \\
C & D
\end{array} \right)
\end{equation}
be a block partitioning of a real matrix $ R$ such that $ A $ is a square $ k \times k $ matrix, $ m \geq n, \; k < n .$ 
Set also
\begin{equation} \label{R_0} 
R_0 = 
\left( \begin{array}{cc}
A & B \\
C & 0
\end{array} \right)
\end{equation}
Clearly 
$
\| R - R_0 \| \leq \|D \| \nonumber
$.
Let $ R_i $ be the optimal rank $ i$ approximation of $ R $ and
let $ R_{0i} $ be the optimal  rank $ i$ approximation of $ R_0.$
According to a well known Weyl inequality (cf. e.g. \cite{Ho}, \cite{Go}), for any matrices $ X, Y $ 
\begin{equation}
\sigma_{i}( X + Y ) - \sigma_{i} ( X ) <= \sigma_1( Y ), \; i = 1,2,\cdots, n  
\end{equation}
Setting here $X \leftarrow R_0 , \; Y \leftarrow D $ and using inequality (3)
twice by interchanging $ R $ and $ R_0 $ we get
\begin{eqnarray}
- \| D \| \leq |\ R - R_i \| - \| R_0 -  R_{0i} \| \leq \| D \|  \nonumber \implies \\
\| R_0 -  R_{0i} \| - \| D \| \leq |\ R - R_i \|  \leq   \| R_0 -  R_{0i} \| + \| D \|
\end{eqnarray}
for $ i = 1, 2,\cdots, n.$ Noticing that
\begin{equation}
\| R - R_{0i} \| = \| R - D - R_{0i} + D \| \leq  \|R_0 - R_{0i} \| + \| D \| \nonumber
\end{equation}
and interchanging here $ R $ and $ R_0 $  we get also  the following  estimate (cf. \cite{Dri}  or \cite{McSherry} Lemma 2)
\begin{equation}
\| R - R_{0i} \| -2 \| D \| \leq \| R - R_{i}  \| \leq  \|R - R_{0i} \| + 2\| D \| 
\end{equation}
In other words, the operator norm error of replacing rank $i$ approximation of $ R $ with rank $i $ approximation of $ R_0 ,$  is no greater than two operator norms of the removed $ m \times (n-k) $ block $ D $ and the estimates (4,5) do not even depend on the fact that $ D $ is  a rectangular matrix block. Further, since rank of $ R_0 $ is no greater than $2k$   it follows from (4) that 
\begin{equation}
\sigma_{i+1} (R) \leq \| D \| \; \; \text{if } i \geq 2k  
\end{equation}
\newline
and it  is quite obvious that 
\begin{eqnarray}
\| R_0 - R_{0k} \| \; \leq \; \min \{\| B \|, \|C\| \}  \\
\| R - R_{k} \| \; \leq \; \min \{\| B \|, \|C\| \} + \| D \| 
\end{eqnarray} 
For an arbitrary matrix $ X $, let  $ Q \equiv Q(X), Q' \equiv Q'(X) $ be its left and right orthogonal SVD multipliers. In other words
\begin{equation} \label{right-rotation}
\; Q R Q' = \Sigma ,\; Q  \equiv Q(R) =  \; \left( \begin{array}{cc}
c_1 & s_1 \\
s_2 & c_2
\end{array} \right)  , \; Q'  \equiv Q'(R) = \left( \begin{array}{cc}
c'_1 & s'_1 \\
s'_2 & c'_2
\end{array} \right)
\end{equation} 
where $\Sigma \equiv \Sigma(R)$ is a diagonal matrix of downward ordered singular values of $ R 
$ and it is assumed that
 matrix blocks $ c_i \equiv c_i(R), s_i \equiv s_i(R),
 c'_i \equiv c'_i(R) , s'_i \equiv s'_i(R), \; i = 1,2 $ are dimension-compatible with the block structure of $ R $ in (1).
Let $(X,Y) $ be a partitioning  of $ \Sigma(R)  =  QRQ'  $ into  $ m \times i $ and $ m \times (n - i) $ vertical matrix bands. Let also 
$ (X_0, Y_0) $ be a similarly shaped partitioning of $ QR_0 Q' .$ 
Clearly
\begin{equation}
\| R_0- R _{0i} \| \leq \| Y_0  \| = \| Y + (Y_0 - Y ) \| \leq \|Y - Y_0 \| + \| Y \|  \
\end{equation}
where $ Y = R - R_i $ by  definition of singular value decomposition. Further,  $ (X- X_0, Y-Y_0 )  \equiv  QRQ' - QR_0Q' \equiv QDQ' $ 
and therefore $ Y - Y_0 $ is just the rightmost $ m  \times ( n-i )$ vertical band of the matrix   
\begin{equation} QDQ' \equiv 
\left( \begin{array}{cc}
s_1 Ds'_2 & s_1Dc^{\prime}_2 \\
 c_2 D s'_2 & c_2Dc^{\prime}_2  
\end{array} \right) = 
\left(\begin{array}{c}
 s_1 \\ c_2 
\end{array} \right) D \;
 ( s'_2, \;  c^{\prime}_2 )  
\end{equation}
It follows then, that  
\begin{equation}
\| R_0 - R_{0i} \| \leq | R - R_i | + \| D (Q'(R)  [k+1:n\;,i+1:n] )\| 
\end{equation}
 and similar estimate is obtained by replacing column partitions with row partitions 
\begin{equation}
\| R_0 - R_{0i} \| \leq | R - R_i | + \| Q(R)[i+1:m,k+1:m]D\| 
\end{equation} 
where "submatrix designation" notation (cf. \cite{Go} p. 27) is used on the right hand sides of (12,13). To get opposite direction inequalities, interchange $ R $ and $ R_0 $ in (13) 
\begin{eqnarray}
 \| R - R_i \| \leq \| R_0 - R_{0i} \| + \| Q(R_0)[i+1:m,k+1:m]D\| \nonumber \\
  | R - R_i | \leq \| R_0 - R_{0i} \| + \| D (Q'(R_0)[k+1:n,i+1:n])\|
\end{eqnarray}
Setting
\begin{eqnarray}
\mu_{i,k}(R) =  \min \{ \; \| D (Q'(R)  [k+1:n\;,i+1:n]) \| , \; \| Q(R)[i+1:m,k+1:m]D \| \;  \}   \nonumber \\
 \bar{\mu}_{i,k}(R) = \max \{ \mu_{i,k}(R), \mu_{i,k}(R_0)\}, \;  \; \; \; \;\;\;\;
\end{eqnarray}
for $i = 1, \cdots n $ and
\begin{eqnarray}
c_2(R_0,i) =   c_2(R_0)[i - k :m-k,1:m-k], \; c'_2(R_0,i) =   c'_2(R_0)[1:n-k,i-k:n-k] \nonumber  \\
c_2(R,i) =   c_2(R)[i-k:k,1:m-k], \; c'_2(R,i) =   c'_2(R)[1:n-k,i-k:n-k] 
\end{eqnarray}
 for $ i > k $   (cf. (9)),
we obtain by the way of introduction 
\begin{theorem} For any $ i = 1,2 , \cdots , n $ 
\begin{equation}
\| \sigma_i(R) - \sigma_i(R_0) \| \leq \ \bar{\mu}_{i,k}(R) 
\end{equation}
In particular,   if $ i > k $ then
\begin{eqnarray}
\| \sigma_i(R) - \sigma_i(R_0) \| \leq \max \{  \min \{ \|c_2(R,i)D\|, \|Dc'_2(R,i)\| \}, 
\min \{ \| c_2(R_0,i)D\|, \|Dc'_2(R_0,i) \|\} \; \} \; \;
\end{eqnarray}
 and if  $ i > 2k $ then 
\begin{eqnarray}
| \sigma_i(R) | \leq \max \{  \min \{ \|c_2(R,i) D\|, \|Dc'_2(R,i)\| \}, 
\min \{ \| c_2(R_0,i) D\|, \| Dc'_2(R_0,i)\} \; \}  \; \; \; \; \;
\end{eqnarray}
\end{theorem}

\begin{example}
	Let 
\begin{equation} \label{right-rotation}
 r  = \left( \begin{array}{cc}
	c & -s \\
	s & c
	\end{array} \right) \nonumber 
\end{equation} 
be a plane rotation matrix and let   
\begin{equation}
r_0  = \left( \begin{array}{cc}
	c & -s \\
	s & 0
	\end{array} \right) \nonumber
\end{equation}
Then $ \sigma_2( r ) = 1 , \; \sigma_2( r_0 ) = 
\left(\frac{ 1 + s^2}{2} - \sqrt{  \left( \frac{1 + s^2}{2} \right)^2 - s^4 }  \right) ^{1/2}, \;  \mu_{1,1}(r) = c $ and  direct computation shows that  $ \sigma_2(r) - \sigma_2( r_0 ) 
= c^2 $ in accordance with Theorem 1.  Computational  experiments show that  $ \max \{ \mu_{i,k}(R) , \mu_{i.k}(R_0) \} $ 
in (15) cannot be replaced with either $ \mu_{i,k}(R) $ or $ \mu_{i,k}(R_0) $.
\end{example}
\begin{example}
Let $ M $ be a symmetric $ p \times p $ matrix and let $ M' $ be the same matrix with the corner element $ M_{p,p} \neq 0 $ zeroed out.  Let $ Q, Q' $ be orthogonal matrices that diagonalize $M $ and $ M' $ respectfully.  Then by Theorem 1 
\begin{equation}
| \sigma_p( M )  - \sigma_p(M') | \; \leq \; | M_{p,p} | \max \{ | Q_{p.p}, |, |Q'_{p,p} |
\}   \nonumber    
\end{equation}
 and therefore, an upper bound equality  $ | \sigma_p( M )  - \sigma_p(M') | = | M_{p,p} | $ is reached if and only if 
$ \sigma_p(M) = | M_{p,p}| ,\; \sigma_p(M') = 0  $.
\end{example}

\noindent In practice, estimates (17,18) are no more useful than (4) since computing numbers $ \mu_{i,k} $ seems to be no easier than diagonalizing original matrix $ R .$ We will explore possibilities of improving estimates (7-8) and (17-19) in the next section. The idea is to block-diagonalize matrix $R$ using block-Givens 
 rotations that can be written down in closed form. A case when the column norms  of the matrix $ R $ are rapidly decreasing is considered in section 3.  
We end the introduction with a few auxiliary lemmas. 
\subsection{Preliminaries }
\begin{lemma}
Let $ f $ be a reasonable non-negative monotone function defined on $ [0, \infty) $ and let 
 $ X $ be a positive semi-definite symmetric $ n \times n $  matrix of rank $ r $
\begin{itemize}
\item[(i)] If $f$ is non-decreasing  then 
$ \sigma_i( f(X)) = f (\sigma_i(X) ) $ for $ i=1,2,\cdots , n$. In particular
$ \| f(X) \| = f( \|X \| ) $
\newline
\item[(ii)] If $ f $ is non-increasing then 
 $ \sigma_i( f(X) ) = f( \sigma_{n-i+1}(X ) ) , \; i = 1,2, \cdots, n $.
 In particular, if $ X $ is of maximal rank ($n=r$)  then                                                                                                                                                                                                                                                                                                                                                                                                                                                                                                                                                                                                                                                                                                                                                                                                                                                                                                                                                                                                                                                                                                                                                                                                                                                                                                                                                                                                                                                                                                                                                                                                                                                                                                                                                                                                                                                                                                                                                                                                                                                                                                                                                                                                                                                                                                                                                                                                                                                                                                                                                                                                                                                                                                                                                                                                                                                                                                                                                                                                                                                                                                                                                                                                                                                                                                                                                                                                                                                                                                                                                                                                                                                                                                                                                                                                                                                                                                                                                                                                                                                                                                                                                                                                                                                                                                                                                                                                                                                                                                                                                                                                                                                                                                                                                                                                                                                                                                                                                                                                                                                                                                                                                                                                                                                                                                                                                                                                                                                                                                                                                                                                                                                                                                                                                                                                                                                                                                                                                                                                                                                                                                                                                                                                                                                                                                                                                                                                                                                                                                                                                                                                                                                                                                                                                                                                                                                                                                                                                                                                                                                                                                                                                                                                                                                                                                                                                                                                                                                                                                                                                                                                                                                                                                                                                                                                                                                                                                                                                                                                                                                                                                                                                                                                                                                                                                                                                                                                                                                                                                                                                                                                                                                                                                                                                                                                                                                                                                                                                                                                                                                                                                                                                                                                                                                                                                                                                                                                                                                                                                                                                                                                                                                                                                                                                                                                                                                                                                                                                                                                                                                                                                                                                                                                                                                                                                                                                                                                                                                                                                                                                                                                                                                                                                                                                                                                                                                                                                                                                                                                                                                                                                                                                                                                                                                                                                                                                                                                                                                                                                                                                                                                                                                                                                                                                                                                                                                                                                                                                                                                                                                                                                                                                                                                                                                                                                                                                                                                                                                                                                                                                                                                                                                                                                                                                                                                                                                                                                                                                                                                                                                                                                                                                                                                                                                                                                                                                                                                                                                                                                                                                                                                                                                                                                                                                                                                                                                                                                                                                                                                                                                                                                                                                                                                                                                                                                                                                                                                                                                                                                                                                                                                                                                                                                                                                                                                                                                                                                                                                                                                                                                                                                                                                                                                                                                                                                                                                                                                                                                                                                                                                                                                                                                                                                                                                                                                                                                                                                                                                                                                                                                                                                                                                                                                                                                                                                                                                                                                                                                                                                                                                                                                                                                                                                                                                                                                                                                                                                                                                                                                                                                                                                                                                                                                                                                                                                                                                                                                                                                                                                                                                                                                                                                                                                                                                                                                                                                                                                                                                                                                                                                                                                                                                                                                                                                                                                                                                                                                                                                                                                                                                                                                                                                                                                                                                                                                                                                                                                                                                                                                                                                                                                                                                                                                                                                                                                                                                                                                                                                                                                                                                                                                                                                                                                                                                                                                                                                                                                                                                                                                                                                                                                                                                                                                                                                                                                                                                                                                                                                                                                                                                                                                                                                                                                                                                                                                                                                                                                                                                                                                                                                                                                                                                                                                                                                                                                                                                                                                                                                                                                                                                                                                                                                                                                                                                                                                                                                                                                                                                                                                                                                                                                                                                                                                                                                                                                                                                                                                                                                                                                                                                                                                                                                                                                                                                                                                                                                                                                                                                                                                                                                                                                                                                                                                                                                                                                                                                                                                                                                                                                                                                                                                                                                                                                                                                                                                                                                                                                                                                                                                                                                                                                                                                                                                                                                                                                                                                                                                                                                                                                                                                                                                                                                                                                                                                                                                                                                                                                                                                                                                                                                                                                                                                                                                                                                                                                                                                                                                                                                                                                                                                                                                                                                                                                                                                                                                                                                                                                                                                                                                                                                                                                                                                                                                                                                                                                                                                                                                                                                                                                                                                                                                                                                                                                                                                                                                                                                                                                                                                                                                                                                                                                                                                                                                                                                                                                                                                                                                                                                                                                                                                                                                                                                                                                                                                                                                                                                                                                                                                                                                                                                                                                                                                                                                                                                                                                                                                                                                                                                                                                                                                                                                                                                                                                                                                                                                                                                                                                                                                                                                                                                                                                                                                                                                                                                                                                                                                                                                                                                                                                                                                                                                                                                                                                                                                                                                                                                                                                                                                                                                                                                                                                                                                                                                                                                                                                                                                                                                                                                                                                                                                                                                                                                                                                                                                                                                                                                                                                                                                                                                                                                                                                                                                                                                                                                                                                                                                                                                                                                                                                                                                                                                                                                                                                                                                                                                                                                                                                                                                                                                                                                                                                                                                                                                                                                                                                                                                                                                                                                                                                                                                                                                                                                                                                                                                                                                                                                                                                                                                                                                                                                                                                                                                                                                                                                                                                                                                                                                                                                                                                                                                                                                                                                                                                                                                                                                                                                                                                                                                                                                                                                                                                                                                                                                                                                                                                                                                                                                                                                                                                                                                                                                                                                                                                                                                                                                                                                                                                                                                                                                                                                                                                                                                                                                                                                                                                                                                                                                                                                                                                                                                                                                                                                                                                                                                                                                                                                                                                                                                                                                                                                                                                                                                                                                                                                                                                                                                                                                                                                                                                                                                                                                                                                                                                                                                                                                                                                                                                                                                                                                                                                                                                                                                                                                                                                                                                                                                                                                                                                                                                                                                                                                                                                                                                                                                                                                                                                                                                                                                                                                                                                                                                                                                                                                                                                                                                                                                                                                                                                                                                                                                                                                                                                                                                                                                                                                                                                                                                                                                                                                                                                                                                                                                                                                                                                                                                                                                                                                                                                                                                                                                                                                                                                                                                                                                                                                                                                                                                                                                                                                                                                                                                                                                                                                                                                                                                                                                                                                                                                                                                                                                                                                                                                                                                                                                                                                                                                                                                                                                                                                                                                                                                                                                                                                                                                                                                                                                                                                                                                                                                                                                                                                                                                                                                                                                                                                                                                                                                                                                                                                                                                                                                                                                                                                                                                                                                                                                                                                                                                                                                                                                                                                                                                                                                                                                                                                                                                                                                                                                                                                                                                                                                                                                                                                                                                                                                                                                                                                                                                                                                                                                                                                                                                                                                                                                                                                                                                                                                                                                                                                                                                                                                                                                                                                                                                                                                                                                                                                                                                                                                                                                                                                                                                                                                                                                                                                                                                                                                                                                                                                                                                                                                                                                                                                                                                                                                                                                                                                                                                                                                                                                                                                                                                                                                                                                                                                                                                                                                                                                                                                                                                                                                                                                                                                                                                                                                                                                                                                                                                                                                                                                                                                                                                                                                                                                                                                                                                                                                                                                                                                                                                                                                                                                                                                                                                                                                                                                                                                                                                                                                                                                                                                                                                                                                                                                                                                                                                                                                                                                                                                                                                                                                                                                                                                                                                                                                                                                                                                                                                                                                                                                                                                                                                                                                                                                                                                                                                                                                                                                                                                                                                                                                                                                                                                                                                                                                                                                                                                                                                                                                                                                                                                                                                                                                                                                                                                                                                                                                                                                                                                                                                                                                                                                                                                                                                                                                                                                                                                                                                                                                                                                                                                                                                                                                                                                                                                                                                                                                                                                                                                                                                                                                                                                                                                                                                                                                                                                                                                                                                                                                                                                                                                                                                                                                                                                                                                                                                                                                                                                                                                                                                                                                                                                                                                                                                                                                                                                                                                                                                                                                                                                                                                                                                                                                                                                                                                                                                                                                                                                                                                                                                                                                                                                                                                                                                                                                                                                                                                                                                                                                                                                                                                                                                                                                                                                                                                                                                                                                                                                                                                                                                                                                                                                                                                                                                                                                                                                                                                                                                                                                                                                                                                                                                                                                                                                                                                                                                                                                                                                                                                                                                                                                                                                                                                                                                                                                                                                                                                                                                                                                                                                                                                                                                                                                                                                                                                                                                                                                                                                                                                                                                                                                                                                                                                                                                                                                                                                                                                                                                                                                                                                                                                                                                                                                                                                                                                                                                                                                                                                                                                                                                                                                                                                                                                                                                                                                                                                                                                                                                                                                                                                                                                                                                                                                                                                                                                                                                                                                                                                                                                                                                                                                                                                                                                                                                                                                                                                                                                                                                                                                                                                                                                                                                                                                                                                                                                                                                                                                                                                                                                                                                                                                                                                                                                                                                                                                                                                                                                                                                                                                                                                                                                                                                                                                                                                                                                                                                                                                                                                                                                                                                                                                                                                                                                                                                                                                                                                                                                                                                                                                                                                                                                                                                                                                                                                                                                                                                                                                                                                                                                                                                                                                                                                                                                                                                                                                                                                                                                                                                                                                                                                                                                                                                                                                                                                                                                                                                                                                                                                                                                                                                                                                                                                                                                                                                                                                                                                                                                                                                                                                                                                                                                                                                                                                                                                                                                                                                                                                                                                                                                                                                                                                                                                                                                                                                                                                                                                                                                                                                                                                                                                                                                                                                                                                                                                                                                                                                                                                                                                                                                                                                                                                                                                                                                                                                                                                                                                                                                                                                                                                                                                                                                                                                                                                                                                                                                                                                                                                                                                                                                                                                                                                                                                                                                                                                                                                                                                                                                                                                                                                                                                                                                                                                                                                                                                                                                                                                                                                                                                                                                                                                                                                                                                                                                                                                                                                                                                                                                                                                                                                                                                                                                                                                                                                                                                                                                                                                                                                                                                                                                                                                                                                                                                                                                                                                                                                                                                                                                                                                                                                                                                                                                                                                                                                                                                                                                                                                                                                                                                                                                                                                                                                                                                                                                                                                                                                                                                                                                                                                                                                                                                                                                                                                                                                                                                                                                                                                                                                                                                                                                                                                                                                                                                                                                                                                                                                                                                                                                                                                                                                                                                                                                                                                                                                                                                                                                                                                                                                                                                                                                                                                                                                                                                                                                                                                                                                                                                                                                                                                                                                                                                                                                                                                                                                                                                                                                                                                                                                                                                                                                                                                                                                                                                                                                                                                                                                                                                                                                                                                                                                                                                                                                                                                                                                                                                                                                                                                                                                                                                                                                                                                                                                                                                                                                                                                                                                                                                                                                                                                                                                                                                                                                                                                                                                                                                                                                                                                                                                                                                                                                                                                                                                                                                                                                                                                                                                                                                                                                                                                                                                                                                                                                                                                                                                                                                                                                                                                                                                                                                                                                                                                                                                                                                                                                                                                                                                                                                                                                                                                                                                                                                                                                                                                                                                                                                                                                                                                                                                                                                                                                                                                                                                                                                                                                                                                                                                                                                                                                                                                                                                                                                                                                                                                                                                                                                                                                                                                                                                                                                                                                                                                                                                                                                                                                                                                                                                                                                                                                                                                                                                                                                                                                                                                                                                                                                                                                                                                                                                                                                                                                                                                                                                                                                                                                                                                                                                                                                                                                                                                                                                                                                                                                                                                                                                                                                                                                                                                                                                                                                                                                                                                                                                                                                                                                                                                                                                                                                                                                                                                                                                                                                                                                                                                                                                                                                                                                                                                                                                                                                                                                                                                                                                                                                                                                                                                                                                                                                                                                                                                                                                                                                                                                                                                                                                                                                                                                                                                                                                                                                                                                                                                                                                                                                                                                                                                                                                                                                                                                                                                                                                                                                                                                                                                                                                                                                                                                                                                                                                                                                                                                                                                                                                                                                                                                                                                                                                                                                                                                                                                                                                                                                                                                                                                                                                                                                                                                                                                                                                                                                                                                                                                                                                                                                                                                                                                                                                                                                                                                                                                                                                                                                                                                                                                                                                                                                                                                                                                                                                                                                                                                                                                                                                                                                                                                                                                                                                                                                                                                                                                                                                                                                                                                                                                                                                                                                                                                                                                                                                                                                                                                                                                                                                                                                                                                                                                                                                                                                                                                                                                                                                                                                                                                                                                                                                                                                                                                                                                                                                                                                                                                                                                                                                                                                                                                                                                                                                                                                                                                                                                                                                                                                                                                                                                                                                                                                                                                                                                                                                                                                                                                                                                                                                                                                                                                                                                                                                                                                                                                                                                                                                                                                                                                                                                                                                                                                                                                                                                                                                                                                                                                                                                                                                                                                                                                                                                                                                                                                                                                                                                                                                                                                                                                                                                                                                                                                                                                                                                                                                                                                                                                                                                                                                                                                                                                                                                                                                                                                                                                                                                                                                                                                                                                                                                                                                                                                                                                                                                                                                                                                                                                                                                                                                                                                                                                                                                                                                                                                                                                                                                                                                                                                                                                                                                  
  $ \|f(X) \| = f(\sigma_r(X))    $ 
\end{itemize}
\end{lemma}
\begin{lemma} \label{Q}
 Take two matrices  $ X $ and $ Y  $ of dimensions  $ (n-k)  \times k $  and $ (n-k) \times (n-k) $ and let $ Q_1 $ and $ Q_2 $ be 
matrices of dimensions $ k \times ( n-k) $ and $ (n-k ) \times (n-k) $.   
Then the norm of the "horizontally-stacked" $ (n-k) \times (n-k) $ matrix $ (XQ_1, YQ_2 )$  is no greater than $ \sqrt{ \| Q_1Q_1^{T} + Q_2 Q_2^{T} \| } \max \{  \|X\| , \| Y \| \} .$
\end{lemma}
\noindent Proof. let $ \lambda =   \max \{  \|X\| , \| Y \| \} .$
For any $ x \in \reals^{(n-k)} $ we have  
\begin{eqnarray}
<(XQ_1 Q_1^{T} X^T + Y Q_2 Q_2^{T} Y^T)x,x> \; = \;
< Q_1 Q_1^{T} X^Tx, X^Tx > + < Q_2 Q_2^{T} Y^Tx, Y^Tx > \nonumber \\ 
\leq \;\; \lambda^2 <(Q_1 Q_1^{T} + Q_2 Q_2^{T})x, x> \;\; \; \leq \; \;\lambda^2 \| Q_1Q_1^{T} + Q_2 Q_2^{T} \| \nonumber
\end{eqnarray}
\begin{corollary}
If under conditions of Lemma \ref{Q} the matrices $ Q_1 $ and $ Q_2 $ 
are such that $  Q_1^{T} Q_2  = 0 $ then
\begin{equation}
  \| (XQ_1, XQ_2 ) \| \; \leq \;  \sqrt{ \max \{  \| Q_1Q_1^{T} \|,  \| Q_2 Q_2^{T} \|\} } \max \{  \|X\| , \| Y \| \} \nonumber
\end{equation}
\end{corollary}
\noindent Indeed, in this case matrices $  Q_1Q_1^{T} $ and $ Q_2 Q_2^{T} $ commute with each other and can be simultaneously diagonalized. Moreover, the resulting diagonal matrices will have no common non-zero elements. 

\begin{lemma} \label{key}
Under conditions of Lemma \ref{Q},
chose  non-negative  numbers $ j \leq k \leq  i \leq n - k $ and 
let $ P $ be a $ (n-k) \times i $ matrix with its first $j$ columns being the same as the first $j $  columns of  $ X $ and the rest $ i - j $ of its columns being the same as the  columns indexed $ k + j +1, \cdots, i $ of $ Y $.
Then $ \| P \|   \leq \max \{  \|X\| , \| Y \| \} .$
\end{lemma}
\noindent The lemma trivially follows from Corollary 1. Chose  $ Q_1 $ to be a matrix with ones on its main diagonal in the first  $ j $ columns and zeros everywhere else. 
Select  $ Q_2 $ in the same way so that its only nonzero elements are located on the main diagonal in the columns indexed by  $ k+1, \cdots , i    $. 
\newline\newline
For a symmetric matrix $ S $ denote its $ i$-th eigenvalue by $ \lambda_i(S) $. 
\begin{lemma} (see \cite{Fulton} Theorem 2) 
If $  S_1 $ and $ S_2 $ are symmetric $ p \times p $ matrices then $ \lambda_i( S_1  + S_2 ) \; \geq \; \lambda_i( S_1) + 
\lambda_p( S_2) $.  In particular, if 
	 $ S_1 $ and $ S_2 $ are positive semidefinite then $ \sigma_i(  S_1 + S_2 )  
\geq \sigma_i( S_1 ) + \sigma_p( S_2)  $ 
\end{lemma}

\begin{lemma}
For any square $ n \times n $ matrix  $ M $ of rank less than $n$,  there is a  QR-decomposition $ QX = M $ such that the last row of the matrix $ X $ is zero.
\end{lemma}
\noindent Proof.  Housholder QR algorithm (cf. e.g. \cite{Go}) applied to $ M $  will encounter a zero diagonal entry at some point because $ M $ is singular. Thereafter, the remaining Housholder transformations can be chosen in such a way that further diagonal entries will be eliminated. 
\paragraph{}
Below we will use a "block decomposition" of a  block partitioned orthogonal matrix. We will write $ X = X_1 \oplus X_2 \oplus \cdots \oplus X_s $  if all the  non-zero elements of a matrix $ X $ are located in diagonal blocks $ X_1, X_2, ..., X_s .$ A unity $  r \times r $  matrix will be denoted by $ I_r .$  As is well known, any $ n \times n $ orthogonal matrix is orthogonally equivalent to a matrix that looks like $ I_a \oplus -I_b \oplus 
\phi_1 \oplus \phi_2 \oplus \cdots \oplus \phi_{ (n - a - b)/2 } $ where $ \phi_i, i = 1,2, \cdots $ are plane rotations.   
\begin{lemma} (Block-Rotation Decomposition).
	Any ($ n \times n $) orthogonal matrix $ Q $ block partitioned as in (9)  can be factored as 
	\begin{eqnarray} 
	\left(
	\begin{array}{cc}
	q_1  & 0 \\     
	0 & q_2     
	\end{array}    
	\right)  \left[ I_r \oplus  \left(
	\begin{array}{cc}
	c &  -s \\     
	s &  c    
	\end{array}    
	\right)  \oplus I_{n-l} \right] \left(
	\begin{array}{cc}
	q'_1  & 0 \\     
	0 & q'_2     
	\end{array}    
	\right)
	\end{eqnarray}  
	where $ q_1 , q'_1 $ are orthogonal $ k \times k $ matrices, 
	 $ q_2, q'_2 $ are orthogonal $(n-k ) \times (n-k) $ matrices,
	  $ r \leq k, \; l \leq n - k, \; k - r = n-k-l $  and  $ c $ and $ s $ are  $ (k-r) \times (k-r) $ diagonal matrices, $  \| c \| < 1, \; 0 < \|s\| \leq 1 $.
	 Obviously,     
	$ c^2 + s^2 $ must be a unity $ (k-r) \times (k-r) $ matrix
\end{lemma} 
\noindent An almost obvious proof is sketched  here for completeness. Diagonalize blocks  $ s_1, s_2 $ in (9)
using singular value decomposition. This decomposition will almost diagonalize blocks  $ c_1 $ and
$ c_2 $ as well with possible exception of scalar multiples of some square orthogonal sub-blocks along the main diagonal. It is easy to see that these sub-blocks must be symmetric and hence can be diagonalized by block-structure preserving orthogonal transformations as well. Finally, the signs of the  diagonal entries of $ c $ and $s$ can be adjusted by appropriately changing matrices  $ q_i, i = 1,2 $. 
%\begin{remark}
%	Unlike somewhat cumbersome "matrix designation" notation (cf. [3]), block-matrix notation creates some tolerable ambiguity  between a "matrix" and a "matrix block". 
%\end{remark}   
\begin{remark}
	It is easy to see, that for a fixed block structure (9) the  block-rotation decomposition (20)  is unique in the same sense the singular value decomposition is. In particular, absolute values of the entries of  diagonal blocks $ c $ and $ s$ are unique up-to a permutation.    
\end{remark}
\noindent Using factorization (20) and Remark 1 we can define a block partition "weight" $ \omega(Q) \leq 1 $ of a block structure of $ Q $ in (9) as
\begin{equation}
\omega(Q) = \max \{ \| c \|, ||s \| \}  \nonumber
\end{equation}
%\begin{remark}
% Theorem 1 and Lemma 6 that an upper bound ($\|D\|$) in inequality (4) cannot be sharp for $ i \geq n - k $ if  
%\end{remark}
\begin{lemma}
Let $ f $ be a continuous function such that $ f(0 ) =1 $ and let  $ X $ be a symmetric matrix. Let 
$ \nu = \| f(X) |_{\Ima X} \|$. Then  $ \| Y f(X) \|  \leq \nu \|Y\| + ( 1 - \nu ) \| Y |_{\ker X} \|  $ for any (appropriately sized) matrix $ Y $.
\end{lemma}
\noindent Proof. By diagonalizing $ X $ (and hence $f(X)$) we find a symmetric matrix $ M $ and an orthogonal projector matrix $ P$ such that $ MP =PM = 0 $ and $ f(X) = M + P  $. Hence $ \|M \| = \| \nu \| $ and  
 $  \| Yf(X) \|  =  \|  Y ( M + \nu P ) + 
( 1 - \nu )Y P  \|  \; \leq \; \nu \| Y \| + ( 1 - \nu ) \| Y |_{\ker X} \| $.
\newline
\newline
Recall (cf. e.g. \cite{Ho}) that column ($ \| \cdot \|_1$) and row ($ \| \cdot \|_{\infty}$)  matrix norms  are defined as 
\begin{eqnarray}
\max \{\text{column-wise sums of absolute values of matrix entries } \} \nonumber \\ \nonumber
\max \{ \text{row-wise sums of absolute values of matrix entries} \}    \nonumber
\end{eqnarray} 
The following fact that is sometimes called "Schur Test Inequality" (see
\cite{Ho}, \cite{Ma}) is well known
\begin{lemma} 
	$ \|X\| \leq \sqrt{  \| X \|_{\infty} \| X \|_1  } \;  $ for any matrix $ X $. In particular, if $ X $ is symmetric then $ \|X \| \leq \|X\|_{\infty} = \|X\|_1 $.
\end{lemma}

\section{Block-Givens Rotations}
Notational conventions that were introduced in section 1 will be retained till the end of  this section. It is assumed in addition that the block $A $ of the matrix $ R $ in (1) is invertible.  
Let 
\begin{equation} \label{right-rotation}
Q  = \left( \begin{array}{cc}
 c_1 & s_1 \\
s_2 & c_2
\end{array} \right)  
\end{equation}
be an orthogonal $ n \times n $ block matrix with 
$ k \times k $ and $ (n-k) \times (n-k) $ diagonal 
blocks $ c_1, \; c_2 $ and $ k \times (n-k), \; (n-k) \times k $ off-diagonal blocks $ s_1, s_2 .$ 
%We fix the block-structure (8) denoting it as $[k, n-k]$. 
We have
\begin{equation}
s^{T}_1 s_1 + c^{T}_2 c_2 = I_{n-k}  
\end{equation}
and 
\begin{equation}
 R Q \; = \; \left( \begin{array}{cc}
 Ac_1 + Bs_2 & As_1 + Bc_2\\
Cc_1 + Ds_2 & Cs_1 + Dc_2 
\end{array} \right) \nonumber
\end{equation}
It's easy to see that $ Q $ can be chosen in such a way that 
\begin{equation}
As_1 + Bc_2 = 0
\end{equation}    
\begin{lemma}
 If $  As_1 + Bc_2 = 0 $ then $ c_2 $ is invertible.
\end{lemma}
\noindent Indeed, assuming that $ c_2 $ is degenerate, take its QR decomposition $ c_2 = x q $ that satisfies conditions of Lemma 5, so that  $ x $ is a matrix with zero last column.
Multiplying (23) by $ q^T $ on the right, we see that 
the last column of $ As_1 q^T $ must be 
 a zero column, hence last column of   $ s_1 q^T $ is a zero column as well because $A$ is invertible. Therefore the same is true for the last column of the matrix 
\begin{equation}
q_1 = \left( \begin{array}{cc}
 s_1 \\
 c_2
\end{array} \right) q^T  \nonumber
\end{equation}       
and that is, impossible because  $ q_1 $  is orthogonal.
\begin{lemma}  If $  As_1 + Bc_2 = 0 $ then  there is an orthogonal matrix  $ q $ such that
\begin{equation}
c_2 = \left(1 + B^T(AA^T)^{-1}B \right)^{-\frac{1}{2}} \;q   
\end{equation}

%and therefore 
%\begin{equation}
%\| c_2 \| = \sqrt {1 /( 1 + \mu^2 ) } 
%\end{equation}
\end{lemma}
\noindent Proof. Solving (23) for $s_1$ ( $ s_1 = - A^{-1}Bc_2 $ ) we get
\begin{equation}
s^T_1 s_1 = c^T_2(A^{-1}B)^T(A^{-1}B)c_2    \nonumber
\end{equation}
It follows then from (22) that  
\begin{equation}
c^T_2c_2 = I_{n-k} - c^T_2(A^{-1}B)^T(A^{-1}B)c_2  \nonumber
\end{equation}
Using Lemma 9, multiply by $ c^{-1}_2 $ on the right and by $ (c^T_2)^{-1} $ on the left, getting
\begin{equation}
c_2c_2^T = \left( I_{n-k} +  (A^{-1}B)^T(A^{-1}B)  \right)^{-1} \nonumber
\end{equation}
Applying Lemma 9 once again, write polar decomposition of  $ c_2 $ as  $ (c_2c_2^{T})^{\frac{1}{2}} q  $ for some orthogonal $ q $ thus obtaining (24). 
\subsection{Block-partitioned Orthogonal Matrices}
With a little bit of extra work we can write down a generalized Givens rotation matrix of the shape (21) in almost  closed form. Let  $ b_0 q_B = A^{-1}B $ be a polar decomposition of $  A^{-1}B $ so
that $ b_0 $ is a $ k \times k  $  positive semi-definite matrix and $ q_B $
is an orthogonal $ (n - k) \times (n -k ) $ matrix such that  
$ q_B^{T} b_0^{2} q_B =  B^T(AA^T)^{-1}B .$ 
Similarly, let $ q_{C}^T c_0  = CA^{-1} $ be a polar decomposition of $ CA^{-1}.$ The following corollary is obtained by straightforward calculations. 
\begin{corollary} \it{ 
The matrix 
\noindent

\begin{equation}\label{GivMatrix}
\left( \begin{array}{cc}
(1+b_0^{2})^{-\frac{1}{2}} & - A^{-1}B\left(1 + B^T(AA^T)^{-1}B \right)^{-\frac{1}{2}} \\
\\
\left(1 + B^T(AA^T)^{-1}B \right)^{-\frac{1}{2}} (A^{-1}B)^T & \left(1 + B^T(AA^T)^{-1}B \right)^{-\frac{1}{2}} 
\end{array} \right)  \nonumber
\end{equation}  

}
\noindent
is orthogonal  and when applied on the right, annihilates top right corner of $ R .$
\end{corollary}
 The matrix in Corollary 2 is a direct analogue of a two-dimensional Givens rotation
(cf. e.g \cite{Go}). Since it depends on $ A^{-1}B $  only, we can 
set by definition 
\begin{eqnarray} \label{cos}
 b_0 q_B = A^{-1}B    \nonumber \\
\cos(A,B) \; := \; (1+b_0^{2})^{-\frac{1}{2}} \nonumber \\
 \cos(B,A) \; := \; \left(1 + B^T(AA^T)^{-1}B \right)^{-\frac{1}{2}} \equiv q_B^{T}\; (1+b_0^{2})^{-\frac{1}{2}} \; q_B \nonumber \\
 \sin(A,B)  \; :=  A^{-1}B\left(1 + B^T(AA^T)^{-1}B \right)^{-\frac{1}{2}} \equiv 
 b_0 (1+b_0^{2})^{-\frac{1}{2}} \; q_B    \nonumber
\end{eqnarray}
\paragraph{} 
and similarly 

\begin{eqnarray} 
 q_{C}^T c_0  = CA^{-1}  \nonumber \\
\cos(A,C) \; := \; (1+c_0^{2})^{-\frac{1}{2}} \nonumber \\
  \cos(C,A) \; := \; \left(1 + C(A^TA)^{-1}C^T \right)^{-\frac{1}{2}}  
\equiv   q_C^{T} \;(1+c_0^{2})^{-\frac{1}{2}} \; q_C
  \nonumber \\
 \sin(C,A)  \; :=  \; (CA^{-1})^T \left(1 + C(A^TA)^{-1}C^T \right)^{-\frac{1}{2}}
 \equiv  c_0 (1+c_0^{2})^{-\frac{1}{2}}\; q_C \nonumber
\end{eqnarray}
It is easy to check that
\begin{eqnarray} \label{givens-eq}
- A\sin(A,B) + B\cos(B,A) = 0   \nonumber  \\
 -\sin(C,A)^TA + \cos(C,A) C = 0   
\end{eqnarray} 
and therefore, we have 
\begin{corollary} 
\it{ Matrices
\begin{eqnarray} \label{GR}
G_R := \left(
\begin{array}{cc}
\cos(A,B) & -\sin(A,B)    \nonumber \\  \\   
\sin(A,B)^T & \cos(B,A)    \nonumber
\end{array}    
\right)  \\ \nonumber \\ \nonumber
G_L := \left(\begin{array}{cc}
\cos(A,C)  & \sin(C,A) \\ \\  
-\sin(C,A)^T & \cos(C,A)  
\end{array} \right)  \nonumber
\end{eqnarray} 
are orthogonal and satisfy  equations (\ref{givens-eq})}. More precisely, there are orthogonal $ k \times k  $ matrices $ q1,q2 $ such that
\begin{eqnarray}
G_L R = \left(\begin{array}{cc}
q_1\sqrt{ A^TA + C^TC} & \cos(A,C) B + \sin(C,A)D  \nonumber \\  \nonumber
\\ \nonumber
0 &   -\sin(C,A)^T B + \cos(C,A) D
\end{array} \right)    \nonumber \\ \nonumber \\ \nonumber \\
 R  G_R = \left(\begin{array}{cc}
\sqrt{ AA^T + CC^T}q_2 & 0 \nonumber \\ \\ \nonumber
C \cos(A,B) + D\sin(A,B)^T  &  - C\sin(A,B) + D\cos(B,A) 
\end{array} \right)    \nonumber
\end{eqnarray}  
\end{corollary}
\noindent For block-rotation decompositions (cf. (20)) of block-Givens rotations  $ G_R $ and $ G_L $ one has immediately
\begin{corollary} If $ A $ is invertible and $ B \neq 0, \; C \neq 0 $ then  
	\begin{eqnarray}
	\omega(G_R) = \max \left\{\frac{1}{\sqrt{1 + \sigma_r(A^{-1}B)^2 } }, 
	\frac{\sigma_1(A^{-1}B)}{\sqrt{1 + \sigma_1(A^{-1}B)^2 } }\right\}    \; <  \; 1 \nonumber \\ \nonumber
	\omega(G_L) = \max \left\{\frac{1}{\sqrt{1 + \sigma_l(CA^{-1})^2 } }, 
	\frac{\sigma_1(CA^{-1})}{\sqrt{1 + \sigma_1(CA^{-1})^2 } }\right\} \; < \; 1 \nonumber
	\end{eqnarray} 
where $ r = \text{rank}(B), \; l = \text{rank}( C ) $. 
\end{corollary}

\begin{example}
%	The left rotation matrix $ G_L $ in (29) depends only on the left column 
%	$\left(\begin{array}{c}
%	A \\ C 
%	\end{array} \right) $ or $ R $. Therefore, 
	A Housholder transformation  (cf. e.g. \cite{Go}) can be viewed as a particular case of a block-Givens rotation.
Let $ v $ be a column vector and let $ a \in \reals, \; a \neq 0 $. It follows from Corollary 3 that a closed form for a Housholder matrix that "rotates"  vector $\left(\begin{array}{c}a \\ v \end{array} \right) $ into 
$ ( \sqrt{a^2  +  \| v \|^2}, \; 0 , \; \cdots, \; 0 )^T $ 
is
\begin{equation}
\left( \begin{array}{cc}
\left(1+v^Tv/a^2\right)^{-\frac{1}{2}} &  (v/a)^T \left(1 + (v/a) (v/a)^T \right)^{-\frac{1}{2}} \\
\\
  -\left(1 + (v/a) (v/a)^T \right)^{-\frac{1}{2}} (v/a)   & \ \left(1 + (v/a) (v/a)^T \right)^{-\frac{1}{2}}   \nonumber
\end{array} \right)  
\end{equation}  
To verify this directly, note that quadratic form can be evaluated in any orthogonal coordinates, and in this case one can choose coordinate system in which $ v =  (||v\|, \; 0, \; \cdots, \;0 )^T $  
\end{example}

\noindent Set
\begin{eqnarray}
\nu_1 = \frac{1}{\sqrt{1 + \sigma_r(A^{-1}B)^2 } }  \nonumber  \\ \nonumber 
\nu_2 = \frac{1}{\sqrt{1 + \sigma_l(CA^{-1})^2 } }  \\ \nonumber 
\rho_{1}  =  \nu_{1} \|D\| + ( 1 - \nu_{1} ) \| D|_{\ker B} \|    \\ \nonumber 
\rho_{2}  =  \nu_{2} \|D\| + ( 1 - \nu_{2} ) \| D^T|_{\ker C^T} \|    \nonumber 
\end{eqnarray}
Applying block-Givens rotations of Corollary 3 to  matrices $ R_0, \; R $  and using lemmas 1 and 7 we improve estimates  (7-8) as follows 
\begin{theorem} If matrix $A$ in (1) is invertible, then
\begin{eqnarray}
\sigma_{k+1}(R_0)   \leq  
\min \left\{ \frac{\|CA^{-1}\|} {\sqrt{1 + \|CA^{-1}\|^2  } } \| B \|,
\; \frac{\|A^{-1}B\|}{\sqrt{1 + \| A^{-1}B\|^2 } } \| C \| \right\} \leq  
\frac{\|B\| \|C \|}{ \sqrt{ ( \sigma_k(A)^2 + \max\{ \|B\|, \|C\| \}^2 } }   \nonumber \\ \nonumber
\sigma_{k+1}(R)   \leq  
\min \left\{ \frac{\|CA^{-1}\|} {\sqrt{1 + \|CA^{-1}\|^2  } } \| B \| +  \rho_2,
\; \frac{\|A^{-1}B\|}{\sqrt{1 + \| A^{-1}B\|^2 } } \| C \| + \rho_1 \right\}  \\  \nonumber 
\end{eqnarray}
\end{theorem}
\begin{corollary}
	If $ n = m = 2k $ and matrices $A, \; B$ and $C $ are invertible, then
\begin{equation}
\sigma_{k+1}(R)  \leq 
\min \left\{\frac{\|B\| \|C \|}{ \sqrt{ ( \sigma_k(A)^2 +  \|B\|^2 } }  +  \frac{\|A\| \|D\|}{ \sqrt{ ( \|A\|^2 + \sigma_k(B)^2  } } ,
\; \frac{\|B\| \|C \|}{ \sqrt{ ( \sigma_k(A)^2 +  \|C\|^2 } }  +  \frac{\|A\| \|D\|}{ \sqrt{ ( \|A\|^2 + \sigma_k(C)^2  } } \right\}  \\  \nonumber 
\end{equation}
\end{corollary}

\subsection{Block Diagonalization }
Another straightforward application of Corollary 3 is a block-diagonalization algorithm

		\begin{algo} 	A matrix $ R ,$ block-partitioned as in (1) can be block diagonalized by a following iterative procedure.
			\begin{itemize}
				\item[1)] Eliminate bottom left $C$-block of $ R $ by computing   $ R_1 = G_L R $   
				\item[2)] Compute matrix $ R_2 $ by applying right block-Givens rotation that  eliminates top-right $B$-block of $ R_1 $
				\item[3)] Repeat the above steps for $ t = 2, 3, \cdots $ alternatively multiplying $ R_t$ by a left block-Givens rotation when $ t $ is even and multiplying $ R_t$ by a right  block-Givens rotation when $ t $ is odd  
				  
			\end{itemize} 
		\end{algo}
	
\noindent By design 
\begin{equation} \label{right-rotation}
R_t =  \left( \begin{array}{cc}
A_t & B_t \\
0  &  D_t 
\end{array} \right)  , \;  \text{for odd } t \text{ and } R_t =  \left( \begin{array}{cc}
A_t & 0 \\
C_t & D_t
\end{array} \right) \text{for even} \; t, \; t = 1,2,\cdots      \nonumber
\end{equation} 
Next lemma  summarizes some useful properties of the block-diagonalization algorithm.
\begin{lemma} \
	\begin{itemize}
		\item[(i)]  $ \sigma_i(A_{t+1}) \geq  \sigma_i( A_t) $ for all $ i = 1,2, \cdots, k $ and $ t = 0,1, \cdots $ 
		\item[(ii)] $ C_1 = 0 ,\; \sigma_i( A_1 ) = \sigma_i( R[1:m,1:k] ) , 
\; i = 1,2, \cdots k  $
		\item[(iii)] $ \| R[1:m,k+1:n] \|= \| R_1[1:m,k+1:n] \|$ 
		\item[(iv)] For all $ t = 1,2, \cdots $ 
\begin{eqnarray}
	\| D_{t+1} \| \; \leq \;  ( 1 + \|C_t \|^2 / \|A_t \|^2 ) ^{-1/2} \|D_t\| \; \text{if t is even} \nonumber \\
 \| D_{t+1} \| \; \leq \; ( 1 + \|B_t \|^2 / \|A_t \|^2 ) ^{-1/2} \|D_t\| \; \text{if t is odd } \nonumber  \\
 	\| B_{t+1} \| \; \leq \;   ( \sigma_k(A_t)^2 + \|C_t \|^2  ) ^{-1/2} \|C_t\| \| D_t \|  \; \text{if t is even}  \nonumber\\
 	\| C_{t+1} \| \; \leq \; ( \sigma_k(A_t)^2 + \|B_t \|^2  ) ^{-1/2} \|B_t\| \|D_t \| \; \text{if t is odd } \nonumber
	\end{eqnarray}
\item[(v)] If $ \sigma_i( R[1:m,1:k] ) \geq \| R[1:m,k+1:n] \|   $ for some $ 1  \leq i \leq k $ then \\ $ \sigma_i( R_t[1:m,1:k] ) \geq \| R_t[1:m,k+1:n] \|   $  for all $   t = 1,2,\cdots$
\item[(vi)] The algorithm converges to a block diagonal matrix
\end{itemize}
\end{lemma}
\noindent Proof. Statement (i) follows from Lemma 4, statements (ii) and (iii) and (vi) are obvious, statement (iv) can be verified by a direct computation and  (v) follows from (iii), (iv) and Lemma 4.
\begin{example}
	It is well known that for any $ p \times q $ matrix
	$ Y $ with columns $  v_1, \cdots, v_q $ and for any $ i , \; 1 \; \leq  i \leq  \; q $  
	\begin{eqnarray}
	\sum_{j=1}^i \sigma^2_{j}( Y  ) \;  \geq \;  	\sum_{j=1}^i \|v_j\|^2  \nonumber , \; \;
	\sum_{j=i+1}^q \sigma^2_{j}( Y  )  \; \leq \;  	\sum_{j=i+1}^q \|v_j\|^2  \nonumber
	\end{eqnarray}	
and one of the ways to see that  is to use Algorithm 1 
\end{example}
\noindent 

\noindent As a direct corollary of Lemma 11  we have 
\begin{corollary}
If $ \sigma_i( R[1:m,1:k] ) \geq \| R[1:m,k+1:n] \|   $ for some $ 1  \leq i \leq k $ then  first $ i $ singular values of   the matrix $\lim_{t} A_t $ coincide with first $ i $ singular values of $ R $
\end{corollary}	 
\begin{remark} Suppose that block-diagonalization performed by Algorithm 1 is followed by a singular value decomposition of remaining diagonal blocks. Using the same notation as in Theorem 1 for the resulting orthogonal multipliers and their $ "c, s" $-blocks 
	we can state
\begin{theorem1}
	For $ i > k $ there are horizontal $ ( m - i + 1 ) \times ( m - k ) $ slices  $ S, \; S_0  $
	of $ \left(\begin{array}{c}
	s_1(R) \\ c_2(R) 
	\end{array} \right) , \;  \left(\begin{array}{c}
	s_1(R_0) \\ c_2(R_0) 
	\end{array} \right) $  
	and vertical $ ( n -k ) \times ( n - i + 1 ) $ slices $ S', S'_0 $ of 
	$  ( s'_2(R), \;  c'_2(R) ),  \;   ( s'_2(R_0), \;  c'_2(R_0) )  $  such that         
	\begin{equation}
	| \sigma_i(R) - \sigma_i(R_0) | \leq \max \{  \min \{ \| S D\|, \|D S'\| \}, \;
	\min \{ \| S_0D\|, \|DS'_0 \|\} \; \}  \nonumber
	\end{equation} 
\end{theorem1}
\noindent Sketch of the proof. Repeat  the proof of  Theorem 1 and use Lemma 3.    
\end{remark}

\section{ Singular Values of Sparse non-negative Random Matrices }
A matrix  (or a vector) will be  called non-negative if all its entries are non-negative. We will  consider large sparse non-negative random matrices with not too different non-zero entries. The notions  involved will be introduced step by step.  
\paragraph{Matrix Density.} Let $ X $ be $ m \times k $ non-negative matrix. 
For any row index subset $ I' \subset   \{1,2, \cdots , m \} $ and any column index subset 
$ J' \subset \{1,2, \cdots , k \} $  the density of $ |I'| \times |J'| $ sub-matrix $ X[I',J'] $ that is formed by intersection of rows indexed by $ I' $ and columns indexed by $J'$ can be  defined 
 (cf. \cite{Al}) as  
\begin{equation}
 \delta_{I',J'} = \delta(X[I',J']) = 
 \frac{ \sum_{i \in I', j \in J' } X_{ij} } { |I'|  |J' | }   
\end{equation} 
  We define the size $ |x| $ of a non-negative vector $   x = (x_1, \cdots, x_m  ) $,   as a sum of  its coordinates. 
For $ I \equiv I_m = \{1,2, \cdots, m \} $ and  $ J \equiv J_k =  \{1,2, \cdots , k   \}$ we 
will use a shorthand $ \delta = \delta( X[I_m, J_k]) = \delta_{ I_m, J_k} (X)  =  \delta (X)$ 
%In other words, $ \delta(X) $ is defined as the density of the first $ k$ columns of $ X .$ 
\paragraph{Invariant Random Vectors.} 
Let $ G $ be a compact subgroup of orthogonal group $ O(p) $. A random vector $ x $ in $\reals^p $  will be called $G$-invariant if for any $ v \in \reals^p  $  and $ g \in G,  \;\;  E (<gv,x> ) = E(<v,x>) $.  For example a uniformly distributed (on a sphere) random vector is $ O(p) $ invariant.
\begin{lemma} 
Let $V =\reals^p $ be a real $p$-dimensional vector space with the standard  Eucledean norm and let $ G $ be a compact (e.g. finite)  Lie-subgroup of the group of isometries $ O(V) $. Suppose that natural representation of $ G $ in $ V $ does not have fixed points. 
%Let $x, y $ be independent random $G$-equivariant (i.e $G-$invariantly distributed)  vectors in $ V $.
Then
\begin{itemize}
	\item[(i)]  If $ z $ is a  $ G$-invariant random vector in V then $ E ( \mathcal{L}(z)  ) = 0 $ for any linear function $ \mathcal{L} $ on $ V $ 
	\item[(ii)]   $ E (<x,y> )  = 0 $ for any pair of independent random $G$-invariant    vectors  $x, y $ in $ V $ 
	\item[(iii)]  (cf. \cite{Al2}) If in addition $G$-action on $ V$ is absolutely irreducible, then for  independent  $G$-invariant  random  vectors  $x, y $ on a unit sphere in $ V $  we have
	%and $G $ contains the full symmetric group $\Gamma_p $ of %permutations of coordinates.   
	\begin{equation}
E\left(<x,y>^2\right) = \frac{1}{p}  \nonumber
	\end{equation} 
\end{itemize}
\end{lemma}
\noindent Proof. 
The expectation vector $ E(z) = ( E(z_1), \cdots, E(z_n) )  $ must be $ G$-invariant and (i) follows from the linearity of expectation: $ E(\mathcal{L}(z) ) = \mathcal{L}(E(z)) = 0 $. Moreover, by independence of $ x $ and $ y $, $ E(<x,y>) = \;< E(x),E(y) > \;= 0 $ (e.g. by (i)).
\paragraph{}

We start the proof of (iii) by noticing that conditions and conclusions of the lemma do no depend on the choice of coordinates in $ V $ . Take an orthogonal  basis $ e_1,\cdots , e_p \in V $.
An orbit of $ e_1 $ contains exactly $ p $ linearly independent vectors $ f_1 = e_1 = g_1e_1,  f_2 = g_2 e_1, \cdots ,  f_p = g_p e_1 $ for some $ g_1 = 1, g_2, \cdots , g_p \in G $ because $G$-action on $ V $ is irreducible. Applying polar decomposition to a non-degenerate matrix    
$ (\; <e_i, f_j >, \; i, \; j  = 2, \cdots, p \; )  $  if necessary,  we can choose an orthogonal basis  $ e_1, e'_2, \cdots e'_n $  in such a way that $ (p-1) \times (p-1) $   matrix $  M = (\;<e'_i, f_j >, \; i, \; j  = 2, \cdots, p \;)  $ is symmetric positive definite and therefore by Schur product theorem  (see e.g \cite{Ho}), its Hadamard square $ M \circ M =  (\;<e'_i, f_j >^2, \; i, \; j  = 2, \cdots, p \;) $ is non-degenerate. Therefore, we can  fix the orthogonal basis 
$ e_1,\cdots , e_p \in V $ in such a way that the matrix
$ (\;<e_i, f_j >^2, \; i, \; j  = 2, \cdots, p \;) $ has rank $ p - 1 $. 
\paragraph{}
Let  
$ x = \sum_{i}^{p} x_i e_i , \; y = \sum_{i}^{p} y_i e_i$ so that
\begin{equation}
 <x,y >^2 \; =  \sum_{i=1}^p x^2_{i} y^2_{i} \;  +  \; 2 <  \sum_{i \neq j}^p x_{i} x_{j}, \; \sum_{i \neq j}^p y_{i} y_{j} >
\end{equation}

\noindent Since action of $ G $ on $ V $ is absolutely irreducible, the symmetric  square $  S^2(V) $ of $ V $ splits into an orthogonal sum 
\begin{equation}
S^2(V) =  V' \oplus W \nonumber
\end{equation}
of one-dimensional trivial $G$-representation  $ V' \equiv \reals \cdot \sum_i e_i \otimes e_i  $ and a fixed point free  $G$-representation $ W $ (cf. e.g. \cite{Wo}).
 It is easy to see that $\frac{p(p-1)}{2}$-dimensional vector $( x_i x_j , i \neq j) $ can be viewed as a random 
$G$-invariant vector on $ W $ and therefore it follows from (ii) that  expectation of the second summand in (38) is zero.
As was already mentioned, the $ G$-orbit of $ e_1 $ contains exactly $ p $ linearly independent vectors $ f_1 = e_1 = g_1e_1, \; f_2 = g_2 e_1, \; \cdots , \; f_p = g_p e_1 $, for some $ g_1 = 1, g_2, \cdots , g_p \in G $, and  
 \begin{eqnarray}
  E (<x,\; e_1 >^2)  \; \equiv \; E ( <x, \; f_1 > ^2 ) \; = \; E( <x, f_i>^2),  \;\; i = 2, \cdots ,  p 
 \end{eqnarray} 	
 because random vector $ x $ is $ G$-invariant.  Opening brackets in (28) we get 
 \begin{eqnarray}
 E( <x, f_j>^2) \; \equiv  \;  E( <\sum_{i=1}^p x_i e_i, f_j>^2) \; = \;   \sum_{i=1}^p< e_i, f_j>^2  E(x^2_{i})  \nonumber \\ \; + \;    \sum_{i \neq k}^p< e_i, f_j><e_k, f_j>  E(x_ix_k), \;\; j = 1 \;, \cdots, \;p 
 \end{eqnarray}
 The second term on the right hand side of (29) is equal to zero (by (i)) and hence    $ E(x_{i}^2), \; i = 1, \cdots p $  must satisfy  a system of  linear equation 
 \begin{eqnarray}
 \sum_{i=1}^p< e_i, f_j>^2  E(x^2_{i}) \;   = \; \lambda  \; (\; \equiv  E( <x, f_j>^2 \;) ,\; j = 1, \; \cdots, \; p  
 \end{eqnarray}
 for some $ \lambda \in \reals $. In particular,
 since $    \sum_{i=1}^p< e_i, f_j>^2 \; =   1 $, the system of $p$ linear equations 
 \begin{equation}
 \sum_{i=1}^p< e_i, f_j>^2  E(x^2_{i}) \; = \;  \frac{1}{p} , \;\;  j = 1, \; \cdots, \;  p  
 \end{equation}
 has a solution 
 \begin{equation}
 E(x_{i}^2) = 1/p , \; i = 1, \; \cdots, \;  p
 \end{equation} 
As was explained above, the rank of the system (31) is $ p - 1$ and therefore  "the choice of expectations" (32) must be unique.
 It follows   that the expectation of any squared coordinate of $x$ (and of  $ y $)  is equal to  $ 1/p $ and using independence of $ x $ and $ y $ we get (cf. \cite{Al2}):
\begin{equation}
 E(\sum_{i=1}^n x^2_{i} y^2_{i}) = p \cdot \frac{1}{p} \cdot \frac{1}{p} = \frac{1}{p} \nonumber
\end{equation}
 \begin{remark} The main idea behind Lemma 12 comes from  \cite{Al2} where statement (iii)  is established for the case  $ G = O(V) $.
\end{remark}
 \begin{remark} Obviously,  the statement (iii) of lemma 12 can be "norm-scaled" for a pair of fixed norm random vectors $ x, y $ , i. e. $ E\left(<x,y>^2\right) = \frac{1}{p} \|x\|^2 \|y\|^2 $. It is less obvious, but still easy to see that Lemma 12 can be generalized in two ways. First, a range of the random vector $ x $ can be assumed to be  a $G$-invariant  subset of $ V $ (or  a  subset of a (unit) sphere in $ V$) since the  proof presented above depends only on the properties of absolutely irreducible representations and invariance of  random vectors. Second,   the condition $ \| x \| = \|y \| =  1 $ of Lemma 12 (iii) can be "replaced by Fubini's theorem", so that
\begin{equation}
  E\left(<x,y>^2\right) = \frac{1}{p}  E(\|x\|^2) E(\|y\|^2) \nonumber
\end{equation} 
for  random (invariant) vectors with varying norms.
\end{remark} 
\paragraph{Non-negative vectors with fixed norms and  sizes.} 
An "end-point" of a uniformly distributed random nonnegative vector $   x = (x_1, \cdots, x_m  ) \in \reals^m $ with a fixed size $ |x| = a $  belongs to a simplex $ \sigma_x $  that is cut-off from the positive ortant by a hyperplane 
 $ H_x $ defined by an  equation $ \sum_{i=1}^n x_i = |x| \equiv a  $ . Take a normal to $ H_x$ vector $ \eta $ with coordinates $ (1, \cdots , 1 ) $. We have $ E(x) \equiv c_x = (|x|/m) \eta \equiv (a/m)\eta $. Vector $ c_x$ is orthogonal to $ H_x $, so that 
$ x = c_x + r_x $ 
where $ r_x \perp c_x $. Let's fix the length (Eucledean norm) of $x $ as well, requiring that the endpoint of $x $ belongs to a sphere  $ S_x $ given by an equation $ \| x \| = b $, thus confining the range of $ x $ to an intersection 
$ S' =  S_x \bigcap \sigma_x $. 
\begin{remark}
It is easy to see that $S'$ is a sphere if and only if 
%the length of $r_x$ is no greater than the radius of a sphere inscribed into the simplex $\sigma_x $, i. e. if and only if 
\begin{equation}
\|r_x \| \equiv  \sqrt{  \|x\|^2 - \frac{|x|^2}{m} } \; \leq \;  \sqrt{\frac{1}{m(m-1) }} \; = \text{ the radius of a sphere inscribed in } \sigma_x \nonumber
\end{equation}
and in general, the set $ S' \subset \sigma_x $ is non-empty as long as   $ \|x \|  \; \leq \;  a = | x|  $ % \nonumber
\end{remark}
\noindent Take a subgroup $ O_{\eta} $ of orthogonal group $ O(m) $ that leaves the vector $\eta $ fixed.  
  The group $ G $ of symmetries of the simplex ($ \sigma_x $) as a subgroup of   $ O_{\eta}  $  is isomorphic to a full symmetric group of coordinate permutations $ \Gamma_m $.  
 According to the remark above, although  the set $ S' $ is not necessarily   $ O_{\eta} $-invariant, it is nevertheless always  $ \Gamma_m$-invariant and hence restricting (truncating in a sense of multivariate distributions) a random vector $ x $  uniformly distributed on a sphere $ \{ |x | = a, \; \| x\| = b \} $ to the intersection set $ S' $ results in a $ \Gamma_m$-invariant random vector in $ S' $. 
   
%
%
% 
%If $ y $ is a non-negative random vector with endpoint in the  hyperplane $ H_y = \{ \sum_{i=1}^n y_i = |y| \equiv y_0  \}$   then  $ y = c_y + r_y, \; r_y \perp c_y, \; c_y = (|y|/m)\eta $. 
\begin{lemma}  Suppose that random vectors $ x $ and $ y $ are independent and $\Gamma_m$-invariant.  If the sizes and norms of random vectors $ x, \; y $ are fixed then 	
\begin{itemize}	
\item[(a)] 
\begin{eqnarray}
	E( <r_x,r_y>) =  0 \nonumber \\ 
	  E(<x,y>) =  |x||y| / m      
\end{eqnarray}
\item[(b)] 
	\begin{equation}  
	Var( <x,y>) = E (<r_x,r_y>^2) \; = \; \frac{1}{m-1} \nonumber \left( \|x\|^2 - \frac{|x|^2}{m} \right) \left( \|y \|^2 - \frac{|y|^2}{m} \right) \nonumber
	\end{equation}
and therefore by Jenssen's inequality 
\begin{equation}
E( |<r_x,r_y>| ) \leq\frac{1}{\sqrt{m-1}}  \sqrt{ \left( \|x\|^2 - \frac{|x|^2}{m} \right) \left( \|y \|^2 - \frac{|y|^2}{m} \right) } \nonumber 
\end{equation} 
\end{itemize} 
\end{lemma}
\noindent Proof. Statement (a) follows from the orthogonality of $ c_x $ and $r_x $, Lemma 12 (i) and Remark 4.
In more details 
\begin{equation}
<x,y> \; = \; <c_x + r_x, \; c_y + r_y> \;  = \; <c_x, \; c_y > + <r_x,\; r_y >  \nonumber
\end{equation}
and therefore 
\begin{equation}
E ( <x,y> )\; = \; E(<c_x,c_y>) + E(<r_x,r_y> )  = |x||y| / m  \nonumber 
\end{equation}
It is well known that the action of the symmetric group $ \Gamma_{m}  $ on the subspace orthogonal to $ \eta $  is absolutely irreducible. Hence, the  Statement (b) follows from Lemma 12 (ii-iii) and Remark 4:
\begin{eqnarray}
Var(<x,y>) \;  = \;	E( <x,y>^2 ) \; - \; E(<x,y>)^2 \; = \;\;  <c_x,c_y>^2 \;  + \nonumber \\
\; 2 <c_x,c_y> E(<r_x,r_y >)  \; + \;    E(<r_x,r_y>^2 \; \; - \; <c_x,c_y>^2  \;\; = \nonumber \\  E(<r_x,r_y>^2   \; = \; \frac{1}{m-1} \|r_x \| ^2  \|r_y \| ^2 \; = \;
\frac{1}{m-1} \left( \|x\|^2 \;  - \;  \frac{|x|^2}{m} \right) \left( \|y \|^2 \; - \; \frac{|y|^2}{m} \right)   \nonumber
\end{eqnarray}
\noindent\noindent
\begin{remark}\
\begin{itemize}  
\item[(i)] As was already mentioned in Remark 4, Lemma 13 remains generally valid even if the sizes and norms of independent random vectors $ x, \; y $ are not fixed. Sizes and norms on the right hand side of (a) and (b) can be replaced by expectations, for example, assuming only sizes of $ x, y $ are fixed, we get under reasonable conditions  
\begin{equation}E(<r_x,r_y>^2   \; = \;
\frac{1}{m-1} \left( E( \|x\|^2 )\;  - \;  \frac{|x|^2}{m} \right) 
\left( E(\|y \|^2 ) \; - \; \frac{|y|^2}{m} \right)   \nonumber
\end{equation}
\item[(ii)] The statement (a) is essentially borrowed from  \cite{Gold}. A zero-one vector  $ u $  can be viewed  as subset of a set its indexes. If $ v $ is another zero-one vector of the same dimension and if  
$ u $ and $ v $  represent random independently selected  subsets of the coordinate index set  then (33) turns into     (see \cite{Gold})
\begin{equation}
E(<u,v>) = E(| u \cap v |) = \frac{|u| |v|} { m}  \nonumber
\end{equation} 
\end{itemize}
\end{remark}
\paragraph{Random Matrix with fixed column norms and sizes.}   
With Lemma 13 in mind we now  fix sizes $  |u_1|, |u_2|, \cdots , |u_k| $ 
 and \it{norms}   $  \|u_1 \| , \| u_2 \| , \cdots , \|u_k \| $ 
\normalfont of all columns of a non-negative matrix $ X $. In particular,  the matrix density $ \delta = (1/k)\sum_{i=1}^k (|u_k| / m ) $ is thereby also fixed. 
Set
 \begin{eqnarray}
   E = E_k = ( \xi_1, \cdots \xi_k )^T, \;\;   \xi_i = |u_i | / \|u_i\|^2, \; i = 1, \cdots, k  \nonumber \\
 U =  U_k = ( |u_1|, \cdots , |u_k| )  , \; \; D \equiv D_k =   \text{diag}(   \|u_1\|^2  , \cdots,  \| u_k \|^2  )   \nonumber  \\
 H  \equiv H_k = \text{diag}\left( 1  - \frac{|u_1|^2} { m \|u_1\| ^2  },  \cdots, 1 - \frac{ |u_k|^2 } {m \|u_k\| ^2}   \right)  
 \end{eqnarray}
  
\paragraph{Moment Ratio.}  Coefficient of variation of a random variable is defined as a ratio of its standard deviation  to  its 
mean. A ratio  of the square root  of the second moment of a random variable   to its first moment  will be somewhat loosely called below a moment ratio. The moment ratio $ \rho $ is related to the coefficient of variation $ \nu $ by a simple rule $ \rho^2 = \nu^2 + 1 $.  Similarly, we can define a sample moment ratio, so that
 if  $ \psi $ is a sample coefficient of variation of a positive sequence $ a = a_1, \cdots , a_t $  then  
\begin{equation}
\rho(a)  =  \frac{ 	\sqrt{\frac{1}{t} \sum_{i=1}^t a_i^{2} } }{ 	\frac{1}{t}\sum_{i=1}^t a_i }   = \sqrt{1 + \psi^2}  ,
\; \; \psi^2  = \rho^2 - 1      
\end{equation}    
\begin{remark}
 We consider here only positive random variables and/or positive sequences  that significantly differ from $0$ (see below). Such a sequence has  a small sample coefficient of variation  if and only if a ratio of square root of its second moment to its first moment is close to 1. 
\end{remark}
\noindent So far we did not make any assumptions on sparsity or tallness of the matrix $ X$. 
Below we will assume that matrix $ X $ is "sparse enough" and that "size-to square-of-the-norm" ratio vector
$ E_k $ has small coefficients of variation. Set
\begin{itemize}
\item[] $ C = \max_{1 \leq i \leq k } |u_i | $
\item[]  $ L = \max_{1 \leq i \leq k } l_i $ where
$ l_i > 0 $ is a number of non-zero entries in $i$-th column
%\item[] $ \xi_0 =  \max_{1 \leq i \leq k } \xi_i$ 
\end{itemize}
We assume the following set of conditions:  
\begin{itemize}
\item[($S_1$)] 
\begin{itemize}
\item[(i)] $ \|u_i \|^2 \geq |u_i | ,  \; i = 1, \cdots , k $. % $. \|u_i \|^2/l_i \geq 1, \; i = 1, \cdots , k $.
In other words,  as vectors in $ \reals^m $ all the columns of $ X $   are located outside of the interior of a ball of radius $ \sqrt{m}/2 $ centered at a point $ ( 1/2, \cdots , 1/2 ) $ 
% $ l_i > 0 $ is an upper bound  on a number of non-zero entries in $i$-th column.
 % This follows from $ | u_i | \geq l_i  \rightleftarrows |u_i|/l_i \geq 1 $
%\item[(ii)] $ C = \max_{1 \leq i \leq k } |u_i | \leq m$ 
%\item[(iii)]  $ L = \max_{1 \leq i \leq k } l_i  \; < \;  m$
\item[(ii)]  $ C \; \leq \;  m $ 
\item[(iii)] The sample coefficient of variation of $ E_k $ is bounded by
 $ \sqrt{ (1 + \frac{m } {Ck})^2 - 1 } $ (cf. (35) above) or equivalently, the moment ratio of $ E_k $ is bounded by $ 1 + \frac{m}{Ck}$.  
In other words, consider the first moment and the square root of the second moment  of the sequence $ \xi_i, \; i = 1,\cdots, k $  
\begin{equation}
 \Xi_1 \; \equiv \;  \frac{1}{k}\sum_{i=1}^k \xi_i  , \; \; \;  \Xi_2  \; \equiv \;  \sqrt{\frac{1}{k}\sum_{i=1}^k \xi^{2}_i}  \;  \equiv \frac{1}{\sqrt{k} }\|E_k\|  
   \nonumber 
\end{equation}
and require that   
\begin{equation}
 \rho( E_k ) \; = \; \Xi_2/\Xi_1  \;  \leq \; 1 + \frac{m}{Ck}
\end{equation}
\end{itemize}
%  $ \|u_i \|^2/l_i \geq 1  $ where
% $ l_i $ is a number of non-zero entries in $u_i $ and there is a constant $ L << m $ 
\end{itemize}
\begin{remark} \
%The condition (i) implies that $ |u_i | \geq \sqrt{ l_i } \geq 1 $ for all $ i = 1,\cdots, k $. 

\begin{itemize}
	\item[(a)]  All conditions in $(S_1)$ obviously hold for 
	zero-one matrices (with 
	$ C = L $). In particular, for zero-one matrices $ \Xi_1 = \Xi_2  = \rho( E_k  ) = 1 $.
	\item[(b)] A rough estimate of the right hand side of (36)  is $ \frac{m}{Ck} < \frac{1}{k} \frac{1}{\delta} $ 
where $ \delta \equiv \delta_k $ is the density of $ X $
\  
\item[(c)] The condition (i) is satisfied if, for example,  $ \|u_i \|^2/l_i \geq 1, \; i = 1, \cdots , k $.
Indeed,
\begin{equation}
\xi_i = |u_i| / \|u_i\| ^2 = \frac{|u_i|} {l_{i}} 
\frac { \sqrt{l_i} } { \| u_i \|  } 
\frac { \sqrt{l_i} } { \| u_i \|  }   \nonumber
\end{equation}
\item[(d)] It follows from the condition (i) that at least one of the coordinates of $ u_i $ is greater than 1.     
\end{itemize}
		
\end{remark}
Some immediate consequences of conditions  $(S_1)$ are listed below
\begin{lemma} \  
\begin{enumerate} 
	\item[(a)]  $ |u_i | \geq 1, \; \xi_i \leq 1, \; i = 1,\cdots k  $.
	   In particular, $ \| E_k \| \leq \sqrt{k} $  and  \; $ \Xi_1, \; \Xi_2 \; \leq \;  1 $. 
	%$ \Xi_1 \; \leq \;  \xi_0 $
	\item[(b)] $ \frac{|u_i|^2} {m\|u_i\| ^2}
	\leq L/m   , \; i = 1, \cdots, k   $
	\item[(c)] $
\sum_{i=1}^k ( \Xi_2 - \xi_i )|u_i| \leq  m  
$
\end{enumerate} 
\end{lemma}
%\noindent To prove (a), write 
%\begin{equation}
% \xi_i = |u_i| / \|u_i\| ^2 = \frac{|u_i|} {l_{i}} 
%\frac { \sqrt{l_i} } { \| u_i \|  } 
%\frac { \sqrt{l_i} } { \| u_i \|  }   \nonumber
%\end{equation}
%The product of the first two terms is no greater than one due to standard inequality
%between arithmetic and quadratic means.
%The last term is no greater than one by condition $ (i) .$ 
\noindent The statement (a) is obvious (see Remark 8 (d)).
Inequality (b) follows from a standard inequality  between arithmetic and quadratic means 
\begin{equation}
 \frac{|u_i|^2} {m\|u_i\| ^2} = \frac{|u_i|^2} {l_{i}^2} 
\frac {l_i } { \| u_i \|^2  } \frac{l_i}{m}  \nonumber  
\end{equation} 
Using (a) and condition   $(S_1)$  (iii)   we verify  (c) as follows
\begin{equation}
\sum_{i=1}^k ( \Xi_2 - \xi_i )|u_i|  \leq  \sum_{i=1}^k ( \Xi_2 - \xi_i )C  \leq   \frac{ k  \Xi_2 - \sum_{i=1}^k \xi_i }{\sum_{i=1}^k \xi_i } C k   \nonumber
\leq ( 1 + \frac{m}{Ck} - 1 ) Ck = m
\end{equation}
\newline
\newline
Condition $(S_1)$ imposes restrictions on sizes and norms of columns $ ( u_1, \cdots , u_k ) $ 
that constitute the matrix  $ X $. An additional, randomness condition imposed on $ X $ is stated as follows      
\begin{itemize}
	\item[($S_2$)] 
	\begin{itemize}
		\item[(i)] The sizes and norms of the columns  $ ( u_1, \cdots , u_k ) $ of  $ X$ are fixed 
		\begin{equation}
			|u_i | = \mathfrak{s}_i, \;\; \| u_i \| = \mathfrak{n}_i  , \; i = 1,\cdots, k  \nonumber
		\end{equation}
				\item[(ii)] Columns of  $ X $ are non-negative independent random  vectors invariant with respect to the group of permutations of coordinates $ \Gamma_m $ (cf. Lemma 13).
\end{itemize}
\end{itemize}
\begin{remark}
It is easy to see that   Lemma 13 is  valid for random vectors that satisfy conditions $ (S_1)$ and $(S_2)$. The condition $(S_2)$(ii) is satisfied by uniformly distributed random vector, restricted (truncated as a multivariate random variable) to a simplex $ | x | = const $   as was explained above.
\end{remark}

\subsection{Singular Values of the Expected Gram Matrix}
Let now
\begin{equation}
\rho \; \equiv \; \rho(X) \; = \; \frac{ \sqrt{ \frac{1} { k }\sum_{1}^k  \frac{|u_i|^2} { m ^2 } } }     
{  \frac{1}{k}\sum_{i=1}^k \frac{|u_i |} {m}  }
\end{equation}
be a sample moment ratio of the sequence of sizes of the columns of $ X $ and  let $ G = E(X^TX) $ be an expected Gram matrix of $ X $.
\begin{theorem}
	Suppose that random matrix $ X $ satisfies conditions $ (S_1) $ and $(S_2)$.  Let  
	 $\tau  $ be a transposition of the index set $ \{ 1, \cdots, k \}$ that sorts the sequence of norms $ \| u_i \| $ in descending order. Then 
for all $  i = 1,\cdots, k $
\begin{equation}
  ( 1 + k \delta  \rho)   \| u_{\tau(i)} \|^2  \; \geq \;    \sigma_i( G ) \; \geq \;  (1 +  \rho + O(L/m) )^{-1} \| u_{\tau(i)} \|^2  
\end{equation} 
\end{theorem}
\noindent Proof. Let $ M $ be a matrix with zero diagonal and with off-diagonal elements defined by (33), that is let
 $ M_{i,j} = |u_i| |u_j| / m , \; i,j= 1,2, \cdots , k ; \; i \neq j.$  Using condition $ (S_2) $ and definition (34), evaluate the expectation of $ XX^T $ as 
 \begin{equation}
  G = D(  1  +   D^{-1}M ) = D \left( H + E_k \otimes \frac{U_k}{m} \right)   
 \end{equation}
Let 
\begin{equation}
Z =  \left(H + E_k \otimes \frac{U_k}{m} \right) \nonumber
\end{equation}
The following two Lemmas essentially follow from conditions $ (S_1)$.
\begin{lemma}      
\begin{equation}
  \| Z^{-1} \| \; \leq \;  1 \; + \; \rho\; + \; O(L/m)  \nonumber
\end{equation} 
\end{lemma} 
\noindent Proof. By Sherman-Morrison formula (see e.g. \cite{Num})
\begin{eqnarray}
Z^{-1} = H^{-1} - \frac{( H^{-1} E_k \otimes \frac{U_k}{m} H^{-1} ) }
{1 \;  + <\frac{U_k}{m}, H^{-1} E_k > }   \nonumber
%= 
%H^{-1} - \frac{H^{-1} E_k\otimes H^{-1}\frac{U_k}{m} } { 1 + \sum_{i}H_i^{-1}\frac{|u_i|}{m}} \nonumber
\end{eqnarray}
We start by noticing that it follows from Lemma 14 (c) that
\begin{equation}
 \; \Xi_2  \sum_{i=1}^k \frac{|u_i |} {m} \; \leq \; 1 \; + \;  <E_k, U_k/m >  
\end{equation} 
Computing up to the order of $ O(L/m ) $, using condition $(S_1) (ii)$, Lemma 14 (a),(b)  and the fact that the matrices involved are of rank one,  we get

\begin{eqnarray} 
\| Z^{-1} \| \; \leq 1 +  \frac{ \|E_k \| \| U_k/m \| } { 1 + <E_k, U_k/m > } \; + O(L/m) \; \leq \; 
1 +  \frac{ \sqrt{k} \; \| U_k/m \| } { \sum_{i=1}^k |u_i | /m  } \; + O(L/m)\nonumber 
\end{eqnarray}
where the last inequality follows from (40).      
Taking into account the definition of the moment ratio  (37)  we obtain the desired estimate
\begin{eqnarray}
 \| Z^{-1} \| \; \leq \;   
 1 + \frac{ \sqrt{ \frac{1} { k }\sum_{1}^k  \frac{|u_i|^2} { m ^2 } } }     
  {  \frac{1}{k}\sum_{i=1}^k |u_i | /m  } + O( L/m )  
\; = \; 1 + \rho + O(L/m) 
\end{eqnarray} 
%}
\newline\newline
By the same token  we have 
\begin{lemma}  
\begin{equation} 
  \| Z \| \leq   1 + k \delta \rho \nonumber 
\end{equation} 
\end{lemma}
\noindent Indeed, according to (37)  and Lemma 14 (a) 
\begin{eqnarray}
 \| Z \| \; \leq \; 1 +  \frac{  \sqrt{k} \; \|E_k\|} {m}  \sqrt{ \frac{1}{k} \sum_{i=1}^k |u_i|^2 } \; = \;  1 +  \frac{ \sqrt{k} \;\|E_k\|} {m} \; \rho \; \frac{1} {k}   \sum_{i=1}^k |u_i| \;   \leq \nonumber \\   \nonumber  \leq \; 1 + \sqrt{k} \; \| E_k \|  \delta \rho \; \leq \; 1 + k\delta \rho  \nonumber
\end{eqnarray} 
The proof of Thorem 3 is now a simple matter. To verify the second inequality, solve (39) for $ D $, apply Lemma 15 and use the well known 
multiplicative inequality for singular values
$ \sigma_i( X_1 X_2 ) \leq \sigma_i(X_1) \sigma_1(X_2 ) $ (cf. e.g. \cite{Ho}).
The proof of the first inequality is almost the same with Lemma 16 being used instead of Lemma 15.
\begin{remark} \it{ One can say that the expected Gramm matrix $ G $  satisfies a Restricted Isometry Condition in terminology of \cite{Rud}, \cite{Versh2}  } 
\end{remark}

\begin{remark} 
 To get some idea about the order of magnitude of the bounds (38), note that in  most   "practical" cases, $ L $ is much smaller than $m $ and  the moment ratio coefficient $ \rho $ is not significantly larger than $1$ because column sizes of $ X$ are rapidly decreasing.  For example, 
 suppose that $ | u_i | = i, \; i = 1, \cdots , k $. Then $ \rho = 2/\sqrt{3} + O(1/k) $.
 On the other hand, Gershgorin radii of $ G $ can be roughly estimated as $ k\delta |u_i| , \;  i = 1,\cdots k $ where the  average row size $ k \delta $ is usually significantly greater than $1$.      
\end{remark}
\subsection{Some Corollaries}
In view of Example 4, the following statement is not surprising.
\begin{corollary}  In notation of Theorem 3
	\begin{eqnarray}
	\sum_{j=1}^i E(\sigma^2_{j}( X  ))  \; \geq \;    \sum_{j=1}^i \sigma_j( G) \nonumber \\ \nonumber
		\sum_{j=i+1}^k E(\sigma_{j}( X  ))^2 \;  \leq \; ( 1 + \rho + O(L/m))   \sum_{j=i+1}^k \sigma_j( G)  \nonumber
	\end{eqnarray}  
	for all $ i =1,\cdots, k $. 
	%In particular, $ E ( \| X \|^2) \; \geq \;  
%\| G \| $ and   $  E( \sigma_{k}^2(X )) \leq ( 1 + \rho + O(L/m) \sigma_k( G ) $ 
\end{corollary}
\noindent Proof. The first statement does not depend on Theorem 3. Sum of the first $ i $ singular values (known as Ky Fan norm, cf. e.g. \cite{Ho}) is a convex function, hence using Jensen's inequality  we get 
\begin{eqnarray}
\sum_{j=1}^i E(\sigma^2_{j}( X  )) \; = \; E\left(  \sum_{j=1}^i \sigma_{j}( X^TX  ) \right) \; \geq \;  \sum_{j=1}^i \sigma_{j}( E(X^TX)  )   \nonumber 
\end{eqnarray}
On the other hand, by Jensen's inequality and  Theorem 3
(cf. also Example 4)  
\begin{eqnarray}
	\sum_{j=i+1}^k E(\sigma_{j}( X  ))^2 \;  \leq 	\sum_{j=i+1}^k E(\sigma^2_{j}( X  )) \;  \leq  	\sum_{j=i+1}^k \| u_{\tau(j)} \|^2   \; \leq \; ( 1 + \rho + O(L/m) ) 	\sum_{j=i+1}^k \sigma_j( G ) \nonumber
\end{eqnarray}
\newline
\newline
We will use standard notation   $ \| \cdot \|_F $  for the Frobenius norm.
Retaining conditions of Theorem 3 and notations preceding it, set
\begin{equation}
 \mathfrak{N} =  \sum_{p=1}^k \|r_p \| \| X'_{p} \|_F    
\end{equation}
where $ X'_{p} $ is a $ m \times ( k-1) $ matrix with columns $ r_i, \; 1 \leq i \leq k, \; i \neq p  $.
As was mentioned above (cf. Lemma 13), the norms $ \| r_i \| ^2 = \|u_i\|^2 \;  - \;  \frac{|u_i|^2}{m}, \; i = 1. \cdots k  $ are fixed by conditions of Theorem 3.
 
\noindent 
\begin{corollary}  Let 
	$ r_0 = max_{1 \leq i \leq k } \| r_i \| $. For all $ i = 1, \cdots, k $  
\begin{itemize} 
\item[(i)] 
%\begin{equation}
%\sigma_i( G ) - \frac {k-1}{ \sqrt{m-1} }  \; \epsilon^2 \; \leq \; E( \sigma_i^{2}( X) ) \; \leq \;  %\sigma_i( G ) +  \frac {k-1}{ \sqrt{m-1} }  \; \epsilon^2  
%\end{equation}
\begin{equation}
\sigma_i( G ) \; - \; \sqrt{ \frac {k-1}{  m-1}}  \mathfrak{N} \; \leq \;   E( \sigma_i^{2}( X) ) \; \leq    \sigma_i( G ) \;  + \;  \sqrt{ \frac {k-1}{  m-1}} \mathfrak{N}   \nonumber
\end{equation}
\item[(ii)] 
\begin{equation}
\sigma_i( G ) - \frac {k-1}{\sqrt{  m-1}}  \; r_{0}^2   \leq  E( \sigma_i^{2}( X) ) \leq    \sigma_i( G ) +  \frac {k-1}{\sqrt{  m-1}} \; r_{0}^2   \nonumber
\end{equation}
\item[(iii)] (cf. \cite{Versh2}). 
In addition to conditions of Theorem 3, suppose that  $ r_i, \; i = 1, \cdots, \; k  $ are restrictions of uniformly distributed random vectors  (see Remark 9). Then there is a constant $ \mathfrak{ c}  $ that does not depend on $ m, \; k $ and $ X $  such that with high probability  
\begin{equation}
\sigma_i( G ) \; - \; \sqrt{ \frac {k-1}{  m-1}} \; \mathfrak{ c} r_{0}^2  \;  \leq \;   E( \sigma_i^{2}( X) ) \leq    \sigma_i( G ) \; + \;  \sqrt{ \frac {k-1}{  m-1}}\; \mathfrak{ c}  r_{0}^2  \nonumber
\end{equation}
\end{itemize}
\end{corollary}    
\noindent Proof. 

\paragraph{(i)} Write $ X^TX $ = $ G + M $ where $ M_{i,j} = \; <r_i,r_j > $ if $ i \neq j $ and $ M_{i,i} = 0 $ for  $ \; i,j  = 1, \cdots, k $ as in Lemma 13  and let  $  \mathfrak{r}_p =  \sum_{i \neq p}^k |<r_p, \;r_i > |, \; p = 1, \cdots, k $. By (3) and Lemma 8 
\begin{equation}
E( \sigma^2_{i}( X ) ) = E ( \sigma_i( X^TX )  ) \; \leq \; \sigma_i( G ) + E( \| M \| ) \leq E( \max_{p} \;  \mathfrak{r}_p ) \;  \leq \sum_{p} \;  E( \mathfrak{r}_p ) \nonumber
\end{equation}
By Lemma 13 
\begin{equation}
E(\mathfrak{r}_p )  \; \leq \; \| r_p \| \; \frac{ k - 1} {\sqrt { m-1  } } \sqrt{ \frac{1}{k-1}\sum_{i \neq p } \| r_i \|^2 } \;  = \;  \sqrt{ \frac {k-1}{  m-1}} \|r_p \|  \sqrt{ \sum_{i \neq p } \| r_i \|^2 }
\nonumber
\end{equation}
and  (i) follows now for example, from Lemma 4.

\paragraph{(ii)}  Let $ D' $ be a diagonal matrix with diagonal entries 
$ \|r_1 \|, \; \cdots, \; \| r_k \| $ and let  $ M = D' M'  D'  $. As in the proof of (i) (and by Lemma 13)
\begin{equation}
E( \sigma^2_{i}( X ) ) \; \leq \; \sigma_i( G ) +  \| D'\|^2 E( \| M' \| ) \;\leq \;
\sqrt{ \frac {1}{  m-1}}\; E( \max_{p} \;  \mathfrak{r'}_p )  r_{0}  ^2\;  \nonumber
\end{equation}
where now $ \|  \mathfrak{r'}_p  \| \leq k - 1 $.
\paragraph{(iii)} We  follow \cite{Versh2}(Corollary 4.4.8). Let $ T $ be a lower-triangular $ (k-1) \times ( k-1) $ sub-matrix of $ M $, that is let $ T_{i,j} = M_{i,j} $ if  $ i > j $   and $ T_{i,j} = 0 $ if  $ i \leq j $. Write again $ T = D'  T' D'  $. It is easy to see that  matrix entries $ T'_{ij} $ of $ T' $ are independent zero mean random variables with  norms bounded by  $ \sqrt{ \frac{1}{m-1}} $ (Lemma 13). By our assumptions,   matrix elements  $ T'_{ij} $   are distributed as truncated sub-gaussian random variables (see \cite{Versh2}) and hence are themselves sub-gaussian. Therefore, (cf.  \cite{Versh2}, Corollary 4.4.8 for details), we have
\begin{eqnarray}
 \sigma_i( G ) + E( \| M \| ) \;  = \; \sigma_i( G )  +  2 E( T ) \; \leq \; \sigma_i( G )  \; + \; 
\; \sqrt{\frac{k-1}{m-1}}  \; \mathfrak{ c}   r_{0} ^2  \nonumber
% \max_{j}  \sum_{j \neq p}^k |<r_j, \;r_i > | 
\end{eqnarray}

\noindent Finally, combining Corollary 8 with Theorem 3 (iii) we state
\begin{corollary} Under conditions of Corollary 8 (iii) 
\begin{eqnarray}
	( 1 \;  + k \;  \delta  \rho + O(L/m)  )  \; \| u_{\tau(i)} \|^2  \; + \;  \mathfrak{ c} \sqrt{ \frac{k-1}{m-1} }  r_{0}^2   \;  \; \geq \;      E( \sigma_i^{2}( X) ) \; \geq  \nonumber \\ \nonumber
	\; \geq \;  (1 \; +  \; \rho + O(L/m) )^{-1} \;  \| u_{\tau(i)} \|^2 \; - \;
	 \mathfrak{ c} \sqrt{ \frac{k-1}{m-1} }  r_{0}^2  \nonumber
\end{eqnarray}
\end{corollary}

\subsection{ A Generalization}
Let's  relax the condition $(S_2) (i) $ by allowing  column norms of $ X_k$ to vary, retaining, however,  the requirement for column sizes to be fixed. In other word, suppose that the following modification of the condition $(S_2)$ is satisfied
\begin{itemize}
	\item[($S'_2$)] 
	\begin{itemize}
		\item[(i)] The sizes   $ ( u_1, \cdots , u_k ) $ of  $ X$ are fixed 
		\begin{equation}
		|u_i | = \mathfrak{s}_i, \; i = 1,\cdots, k  \nonumber
		\end{equation}
		\item[(ii)] Columns of  $ X $ are non-negative independent random  vectors invariant with respect to the group of permutations of coordinates $ \Gamma_m $ (cf. Lemma 13).
	\end{itemize}
\end{itemize}

\noindent This change in condition $ (S_2 )$ leads  naturally to a change in a definition of the numbers $\xi_i $ so that definitions (34) are replaced by
%\begin{equation}
%  a_i \; \leq \;  \sum_{t=1}^{l_i}  X_{ij_t}^2 \equiv ||u_i\|^2 \; \leq b_i \; 
%\end{equation}  
\begin{eqnarray}
\xi_i = |u_i | / E(\|u_i\|^2), \; i = 1, \cdots, k; \; E_k = ( \xi_1, \cdots \xi_k )^T   \nonumber \\ 
  D \equiv D_k =   \text{diag}(   E(\|u_1\|^2 ) , \cdots,  E(\| u_k \|^2 ) )   \nonumber  \\
 H  \equiv H_k = \text{diag}( 1  - \frac{|u_1|^2} { m E(\|u_1\| ^2 ) },  \cdots, 1 - \frac{ |u_k|^2 } {m E(\|u_k\| ^2} ) ) 
\end{eqnarray}
Note that condition $ (S_1) $ does not require any changes, and therefore nothing happens to the density $ \delta $ and the moment ratio $ \rho$ of the sequence of column sizes of  $ X $. Moreover,  it follows from $ (S_1) $ (i)  that  $  |u_i | / E(\|u_i\|^2) \; \leq \; 1, \; i = 1, \cdots, k $. Therefore, with these adjustments,  Lemma 14 remains valid and the expression (39) for the expected Gram matrix $ G $ of $ X $ does not change. Hence, the proof od Theorem 3 can be repeated verbatim, leading to
\begin{theorem3}
	Suppose that random matrix $ X $ satisfies conditions $ (S_1) $ and $(S'_2)$ with $\xi_i, \; i = 1,\cdots , k  $ defined by (43).  Let  
	$\tau  $ be a transposition of the index set $ \{ 1, \cdots, k \}$ that sorts the sequence of squared column norm expectations $ w_i = E( \| u_i \|^2 ) $ in descending order. Then 
	for all $  i = 1,\cdots, k $
	\begin{equation}
	( 1 + k\delta \rho + O(L/m )  \; w_{\tau(i)} \;   \geq   \; \sigma_i( G ) \; \geq \;  (1 +  \rho + O(L/m) )^{-1} \; w_{\tau(i)}  
	\nonumber
	\end{equation}
%	  (with $\xi_i $ defined by (58))  
\end{theorem3}

\subsection{Gamma Distributed Column Sizes}
Finally, we will briefly touch upon  a case when column sizes  $ |u_i|, \; i = 1,2, \cdots, k  $ are sampled from a known distribution for which an estimate of the sample moment ratio $ \rho $ is available. For gamma distribution, the result of this kind  was  obtained in \cite{Hwang}, we state it here in a form convenient for what follows. 
\begin{lemma} (see \cite{Hwang})
Let the sample  $ Y = (y_1, \cdots ,y_k ) $ be drawn from a population with gamma density  
	\begin{equation}
	g(t,\alpha, \beta) \equiv \frac{\beta^{\alpha}}{\Gamma(\alpha) } t^{\alpha-1}e^{-\beta t}  , \;\; \alpha \geq 1 , \; \beta, \; t  > 0 
	\end{equation}
and let $ A_k, \; S^{2}_k  $ and $ \Theta^2 $  be respectfully its sample mean, sample variance and sample second moment. Then 
\begin{equation}
E \left( \frac {S^{2}_k }{ A^2_{k} } \right) \; = \; \frac{ E(S^{2}_k) } {E(A^2_{k} )} \; = \;\alpha^{-1} + 
O\left( \frac{1}{k} \right)   \nonumber
\end{equation}
and therefore, by (35)
\begin{equation}
E \left( \frac {\Theta^{2}_k }{ A^2_{k} } \right) \; = \; \frac{ E(\Theta^{2}_k) } {E(A^2_{k} )} \; = \; 1 + \alpha^{-1} + 
O\left( \frac{1}{k} \right)   \nonumber 
\end{equation}
and  by Jenssen's inequality
\begin{equation}
E( \rho( Y )) \; = \; E \left( \sqrt{\frac {\Theta^2_k }{ A^2_{k} } } \right) \; \leq \;  \sqrt{ \; 1 + \alpha^{-1} }  + 
O\left( \frac{1}{k} \right) \nonumber
\end{equation}
\end{lemma}

%  ===========================================================================================
%  [The estimate of $ \rho $ in the theorem can be larger (upper bound estimate)  on either side because it's division on the right and mult on the left in Theorem 3 and 3'. hence inequality above is good enough.  ]
%
%=================================================================

\noindent A minor difficulty in combining Lemma 17 with Theorem 3$'$ is that to satisfy condition $(S_1) $ we have to consider a case of a left-truncated  gamma distribution $ g_{\geq a}( t, \alpha, \beta ), \; a \geq 1  $,	 defined as 
\begin{eqnarray}
g_{\geq a}( t, \alpha, \beta ) = 0  \; \; \text{if} \; t \in [0, a] ; \;  \;  
 g_{\geq a}( t, \alpha, \beta ) \; = \; \frac{ g(t,\alpha, \beta) }{1 - F(a)}  \; \;  \text{if} \; t \in [a, \infty]   
\end{eqnarray}
where 
\begin{equation}
F(x) = \frac{\beta^{\alpha}}{\Gamma(\alpha) }\int_{0}^{x}  t^{\alpha-1}e^{-\beta t} dt \nonumber
\end{equation}
It is not hard to estimate sample coefficient of variation for a random variable that is distributed according to (45).  
First, note  that a (gamma) distribution is a mixture of its left-truncated and right-truncated distributions
\begin{equation}
g(t,\alpha,\beta) \; = \; F(a) \; g_{\leq a}( t, \alpha, \beta ) \; + \; ( 1-F(a) ) \; g_{\geq a}( t, \alpha, \beta ) 
\end{equation}
where
\begin{eqnarray}
g_{\leq a}( t, \alpha, \beta ) = 0  \; \; \text{if} \; t \in (a, \infty)  \; \; \text{and} \;\;  
g_{\leq a}( t, \alpha, \beta ) \; = \; \frac{ g(t,\alpha, \beta) }{F(a)}  \; \;  \text{if} \; t \in [0,a] 
\end{eqnarray}
\begin{lemma}
	If $ a = 1, \;  \alpha \geq 1 \;$ and $\; \beta  < 1  $ then the moment ratio $ \rho $ of a random variable that follow left-truncated gamma distribution (45) can be estimated as
\begin{equation}
\rho \; = \; \sqrt{1 \; + \; \frac{1}{\alpha} } \; + \; O( \beta^2 )   
\end{equation} 
\end{lemma}
\noindent Indeed, if $ \alpha \geq 1 $, then  
\begin{equation}
F(1 ) \; = \; \frac{\beta^{\alpha}}{\Gamma(\alpha) }\int_{0}^{1}  t^{\alpha-1}e^{-\beta t} dt \; \leq \; 
\beta \int_{0}^{1} e^{-\beta t} dt \; = \; \beta + O( \beta^2 )   
\end{equation}
\noindent and a similar estimate is valid for moments of a random  variable with  a right-truncated
distribution  (47). Now from (46), we have  
\begin{eqnarray}
F(1) \mu_1 \; + \: ( 1 - F(1) ) \mu_2 \; = \; \frac{\alpha} {\beta}, \; \; \;  
F(1) s_1 \; + \: ( 1 - F(1) ) s_2 \; = \; \frac{\alpha + \alpha^2} {\beta }    \nonumber
\end{eqnarray}
where $ \mu_1, s_1, \; \mu_2,  s_2 $ are first two moments of  left and right truncated distributions (45) and (47)  with $ a = 1 $.  Dividing the second equation by the first and taking into account (49) we get (48).
\begin{corollary} Under conditions  (and in notation) of Theorem $3'$, suppose in addition that 
\begin{enumerate}
	\item[(a)] $ \frac{1}{k} \; \leq \;  \frac{L}{m} $
	\item[(b)] The column sizes $ |u_i| , \; i = 1, \cdots , k  $ are  sampled from  left-truncated gamma distribution (45) with $ a = 1, \; \alpha \geq 1 $ and $ \beta \leq 1/\sqrt{k} $   
\end{enumerate} 
Then with high probability
	\begin{eqnarray}
	\left( 1 +  \frac{k}{m\beta} \sqrt{ \alpha (\alpha + 1 ) }\; + \;  O(L/m \right)  w_{\tau(i)}  \;  \geq \;    \sigma_i( G ) \; \geq \;     \left(1 +  \sqrt{1 + 1/\alpha}  \;  + \;  O(L/m) \right)^{-1}  w_{\tau(i)}  \nonumber
	\end{eqnarray} 
\end{corollary}
\noindent Proof. Start by estimating the  matrix density (cf. \cite{Hwang}):
\begin{equation}
E(\delta ) = E \left( \frac{1}{k}\sum_{i}^k \frac{ | u_i|}{m} \right ) \equiv \frac{1}{m}E \left(   \frac{1}{k}\sum_{i} | u_i|\right)  =   \frac{1}{m}\frac{\alpha} {\beta } %+O\left(\frac{1}{k}\right)  
\nonumber \\
\end{equation}
It follows  from Lemma (18) that we can use the estimate of Lemma 17  for the sample moment ratio of the sequence of column sizes when condition (b) is satisfied. Now, using condition (a)  replace the sample moment ratio   $ \rho$ in Lemma 15 (cf. (41))  by the estimate obtained in Lemma 17. The same substitution can be done in Lemma 16, where just obtained estimate for the matrix density can be used as well.

\begin{remark}
The paper  \cite{Hwang}  contains some information on estimation of  confidence intervals for the sample coefficient of variation. These results can be used to specify more precisely what is meant by "high probability" in Corollary 10
\end{remark}

\section{Tall Sparse Matrix with rapidly decreasing column norms   }
We now return to the subject of first two sections. Let $ R_{left} \equiv R[1:m,1:k]  $ be a left $m \times k $ sub-matrix of  $ R $ in (1) and let $ R_{right} \equiv R[1:m,k+1:n]  $ be the right $ m \times ( n-k) $ sub-matrix of $ R $.  In other words, $ R_{left} $  is  a matrix formed by the first $k$ columns $  u_1, \cdots , u_k $ of $ R $.
It follows from Corollary 6 that the Algorithm 1 will produce top $ i, \; i \leq k  $  singular values of  $ R$ as long as $ \sigma_i \left( R_{left} \right) \geq \| R_{right} \| $. Is  there a way to  establish this estimate beforehand,  without computing the singular value decomposition?  
Suppose that the columns of $ R $ are  sorted in descending order of their norms and that 
\begin{itemize}
\item[] the sizes of the columns of $R $ follow a truncated gamma distribution (45) with the shape parameter $ \alpha $.
\item[]  $ R $  is sparse enough in a sense that constant $ L/m $ of the  Corollary 10 is small.
\item[] the ratios $ |u_i| / \|u_i\| ^2, \; i = 1, \cdots, n $  are roughly equal   
\end{itemize}
\noindent Although Corollary 10 is  not applicable to a $ R $,
we nevertheless,  can adopt an estimate 
\begin{equation}
\sigma_i ( R_{left} )  \; \geq \;
( 1 + \sqrt{ 1 + 1/\alpha}  )^{-1/2} \|u_i \| , \; i < k  \nonumber
\end{equation} 
as a heuristic guidance. 
Note that the norm or the matrix $ R_{right}  $  
is bound from above by $   \sqrt{| u_{k+1}  | \| R_{right} \|_{\infty} } \; $ (cf. Lemma (8)).
%In other words,
%\begin{equation}
%\| R_{right} \|  \leq  \sqrt{ | u_{k+1}  | \times \text{ (the maximal size of rows of }  R_{right}) } 
%\end{equation}
Hence we can assume that the Algorithm 1 will find top $ i $ singular values of $ R $ if there is an index $ i < k $ such that
\begin{equation}
  \|u_i \| \geq  \left(  1 + \sqrt{ 1 + 1/\alpha}  \; \right )^{1/2}  \sqrt{ | u_{k+1}  | \times \text{ (the maximal size of rows of }  R_{right}) } \nonumber
\end{equation}
\begin{example} Take
a $ 10^7 \times 10^5 $ sparse (e.g. with density $ \approx 10^{-4}$) non-negative matrix $R$, partitioned as in (1) with
$ k = 10000 .$ It is reasonable to assume that the sizes or the rows are on average a hundred times smaller than the sizes of the columns. Suppose that distribution of sizes of the first $ 10000 $ column norms of $ R $ is exponential ($ 1 + 1/\alpha  = 2 $)  and let $ i = 1000$.
Assuming that  non-zero matrix entries in the first 10000 columns are spread-up evenly, 
 the Algorithm 1 has a chance to recover at least top 1000 singular values of $ R $ if the norm of the 1000-th 
column of $ R $  is no less than $ 0.1554  $ times the size of the 10001-th column
\end{example} 
\noindent Consider now a following low rank approximation procedure.

\begin{algo} \it{
Given a large sparse $ m \times n $ matrix $R$ 
\begin{itemize}
\item[1)] Find  block partitioning (1) with relatively dense
invertible $ k \times k $ matrix $ A $ and sparse  small-norm matrix $ D $
\item[2)] Throw away the bottom right block $ D $ and use Algorithm 1 to find  a low rank approximation of the remaining matrix $ R_0 .$ 
\end{itemize} }
\end{algo}
\noindent We summarize the heuristics behind this algorithm as follows
\begin{remark} \it{
The Algorithm 2 will approximate  a few top singular values  of a tall sparse non-negative matrix  with an approximation error of less than two operator-norms of the removed bottom right block } 
\end{remark}

\section{Acknowledgments}

The author would like to thank Rama Ramakrishnan for useful comments and stimulating discussions.

\end{document}